\documentclass[12pt,linesnumbered,ruled,vlined]{article}
\usepackage{mathptmx}
\usepackage[utf8x]{inputenc}
\usepackage{color}
\usepackage{array}
\usepackage{booktabs}
\usepackage{mathtools}
\usepackage{amsmath}
\usepackage{amsthm}
\usepackage{amssymb}
\usepackage{graphicx}
\usepackage{multirow} 
\usepackage[unicode=true,pdfusetitle,
 bookmarks=true,bookmarksnumbered=false,bookmarksopen=false,
 breaklinks=false,pdfborder={0 0 1},backref=false,colorlinks=true]
 {hyperref}
\hypersetup{
 pdfborderstyle={},pdfborderstyle={},pdfborderstyle={},urlcolor=blue,linkcolor=blue,citecolor=blue}

\makeatletter

\providecommand{\tabularnewline}{\\}

\@ifundefined{date}{}{\date{}}

\usepackage{amsfonts}
\usepackage{dsfont}

\numberwithin{equation}{section}

\allowdisplaybreaks[4]
\DeclareMathAlphabet{\mathcal}{OMS}{cmsy}{m}{n}

\newcommand{\bfB}{\mathbf{B}}
\newcommand{\bbB}{\mathbb{B}}
\newcommand{\bfX}{\mathbf{X}}
\newcommand{\bbX}{\mathbb{X}}
\newcommand{\bfY}{\mathbf{Y}}
\newcommand{\bbY}{\mathbb{Y}}
\newcommand{\bfZ}{\mathbf{Z}}
\newcommand{\bbZ}{\mathbb{Z}}
\newcommand{\rd}{d_{\mathfrak{R}_1}}
\newcommand{\rrd}{d_{\mathfrak{R}_2}}

\newtheorem{theorem}{Theorem}

\newtheorem{lemma}{Lemma}
\newtheorem{remark}{Remark}
\newtheorem{definition}{Definition}
\newtheorem{proposition}[theorem]{Proposition}
\newtheorem{corollary}[theorem]{Corollary}

\makeatother

\usepackage{authblk}
\usepackage{tikz}

\begin{document}

\title{Rough Path Renormalization from Stratonovich to Itô for Fractional Brownian Motion}

\author[1]{Zhongmin Qian\thanks{qianz@maths.ox.ac.uk}}
\author[2,3]{Xingcheng Xu\thanks{xingcheng.xu18@gmail.com}}

\affil[1]{Mathematical Institute, University of Oxford}
\affil[2]{Shanghai Artificial Intelligence Laboratory}
\affil[3]{School of Mathematical Sciences, Peking University}

\date{}

\maketitle

\begin{abstract}
{\normalsize{}

This paper develops an Itô-type fractional pathwise integration theory for fractional Brownian motion with Hurst parameters \( H \in (\frac{1}{3}, \frac{1}{2}] \), using the Lyons' rough path framework. This approach is designed to fill gaps in conventional stochastic calculus models that fail to account for temporal persistence prevalent in dynamic systems such as those found in economics, finance, and engineering. The pathwise-defined method not only meets the zero expectation criterion but also addresses the challenges of integrating non-semimartingale processes, which traditional Itô calculus cannot handle. We apply this theory to fractional Black-Scholes models and high-dimensional fractional Ornstein-Uhlenbeck processes, illustrating the advantages of this approach. Additionally, the paper discusses the generalization of Itô integrals to rough differential equations (RDE) driven by fBM, emphasizing the necessity of integrand-specific adaptations in the Itô rough path lift for stochastic modeling.

}{\normalsize\par}
\end{abstract}

\vspace{2mm}
\hspace{5mm}\textbf{Keywords: }Rough paths, It\^o integration, fractional Brownian motions, fractional 

\hspace{5mm}Black-Scholes model, fractional Ornstein-Uhlenbeck process, renormalization

\vspace{2mm}

\hspace{5mm}\textbf{MSC{[}2010{]}}: 60H05, 60H10; 91B70, 91G80.

\section{Introduction}\label{intro-sec}

Stochastic calculus, underpinned by Brownian motion and semi-martingales, plays a pivotal role in modeling dynamic stochastic environments across fields such as economics, finance, and engineering. A cornerstone of It\^o calculus is the exploitation of the martingale property, with the It\^o integral defined in relation to the complete law of the driving process, rather than specific realizations. However, empirical analyses of various data often reveal dynamic processes characterized by temporal persistence not accounted for by semi-martingale-driven dynamics. This discrepancy is exemplified by rough volatility models and more broadly, rough path models. A prominent approach to modeling memory effects over time involves the use of fractional Brownian motions (fBMs), a class of Gaussian processes introduced by Mandelbrot and van Ness in 1968 (\cite{MvN68}). These fBMs extend standard Brownian motion by capturing dependencies that decay polynomially over time.

The exploration of integration theories with respect to fractional Brownian motion has garnered significant interest among researchers, as evidenced by a wide range of studies (see, for example, \cite{Bend03a,Bend03b,BHOZ08,CQ02,DU99,EvD03,HO03,NO02} and references therein). Notably, Duncan, Hu, and Pasik-Duncan \cite{DHP00} pioneered the use of the Wick product to define fractional stochastic integrals against fBM, termed Wick-It\^o integrals. Subsequently, Hu and Øksendal \cite{HO03} and Elliott and van der Hoek \cite{EvD03} expanded this methodology, developing a fractional white noise calculus applicable to financial modeling. The method proposed by \cite{DHP00} is applicable specifically to the persistent case of fBM, where the Hurst parameter \( H>\frac{1}{2} \), while the model in \cite{EvD03} addresses the anti-persistent case with \( H<\frac{1}{2} \). These developments employ Wick products instead of ordinary multiplication to define stochastic integrals. Additionally, the concept of Wick products has been adapted for analyzing portfolios and self-financing strategies within fractional Black-Scholes markets. Through innovative concepts, the authors demonstrated the non-arbitrage properties and market completeness when replacing the traditional Black-Scholes model \cite{BS73} with a geometric fBM characterized by a Hurst parameter \( H>\frac{1}{2} \). However, the use of Wick products has sparked significant debate due to conceptual divergences from the conventional Black-Scholes model, particularly regarding the economic interpretation of fundamental concepts like portfolio values and self-financing (see, e.g., \cite{BH05}). It is crucial to understand that Wick products, defined merely as the multiplication of two random variables, lack a pathwise computational basis. Consequently, from the realizations \( X(\omega) \) and \( Y(\omega) \) of two random variables \( X \) and \( Y \), it is not possible to compute the Wick product \( X(\omega) \) and \( Y(\omega) \) directly.

Multiple methodologies have been advanced to define stochastic integrals with respect to fBM. One notable pathwise approach, developed by Ciesielski, Kerkyacharian, and Roynette \cite{CKR93}, and furthered by Z\"ahle \cite{Zah98}, leverages the H\"older continuity of fBM sample paths. Specifically, Ciesielski et al. \cite{CKR93} introduced an integration framework using wavelet expansions within Besov-Orlicz spaces, and Z\"ahle \cite{Zah98} extended this by employing fractional calculus alongside a generalization of the integration by parts formula. These theories, however, are applicable predominantly when the Hurst parameter \( H>\frac{1}{2} \). A second method, pioneered by Decreusefond and Üstünel \cite{DU99}, utilizes Malliavin calculus tailored to fBM. This approach has been extensively explored in subsequent research by Alos and Nualart \cite{AN03}, Carmona, Coutin, and Montseny \cite{CCM03}, and Cheridito and Nualart \cite{CN05}.

A third method, integrating via rough path theory, was introduced by Coutin and Qian \cite{CQ02}. This theory, applicable for fBM with a Hurst parameter \( H>\frac{1}{4} \), aligns with Stratonovich's concept of integration. Despite its theoretical robustness, applying Stratonovich's fractional integration to financial models, such as option pricing, reveals a significant drawback: the expected values of Stratonovich integrals are generally non-zero. This characteristic potentially allows for the construction of arbitrage strategies within fractional Black-Scholes markets (see, e.g., \cite{BSV06,Shir98}). Since fBM does not qualify as a semi-martingale when \( H \neq \frac{1}{2} \), it is infeasible to apply Itô's semi-martingale theory directly to fBM. A more viable goal is to develop a stochastic integral where integrations against fBM yield zero-mean outcomes. Such an approach would be more advantageous than Stratonovich's, as it minimizes systemic bias in pricing processes across all possible paths, thereby mitigating arbitrage opportunities.

This paper endeavors to develop an Itô-type fractional pathwise integration theory tailored to fBM with Hurst parameters in the range \( H \in (\frac{1}{3}, \frac{1}{2}] \), employing a rough path approach. Our proposed theory, which is defined pathwise, achieves zero expectation, aligning with our initial objectives. Rough path theory, as articulated by Lyons, provides a robust framework for handling fBMs with \( H < \frac{1}{2} \) (see, e.g., \cite{FH14,FV10,Lyons98,LCL,LQ02}). Within this framework, we develop essential tools characteristic of Itô-type integration, such as connections to Stratonovich fractional integrals, fractional Itô's formula, chain rule applications, and more. While previous studies have examined rough path integration theory for non-geometric paths (see \cite{Gu10,hairer2015,LY2016}), our focus is on specifically advancing this theory for fBM, aiming to enhance its application in areas like fractional Black-Scholes (fBS) models, fractional Ornstein-Uhlenbeck (fOU) processes, and more broadly, rough differential equations (RDEs) driven by fBMs.

As a practical application, we employ our pathwise Itô integration theory on fractional Black-Scholes models and demonstrate that the corresponding fractional Black-Scholes market remains arbitrage-free under a more restrictive set of trading strategies than those typically permissible in a complete market. While our focus on arbitrage issues illustrates one potential use of this theory, there are numerous other areas where it could be applied. Additionally, we extend the application of our theory to high-dimensional fractional Ornstein-Uhlenbeck processes. This allows us to address estimation challenges associated with fOUs, which are instrumental in modeling and forecasting realized volatility (see also \cite{QX2024}).

The structure of the paper is as follows: Section \ref{sec-prelimi} offers a concise review of essential insights related to rough path theory and fractional Brownian motion. Section \ref{RI-sec} forms the crux of our discussion, where we define integration with respect to the Itô-fractional Brownian rough path, clarify the relationship between Stratonovich fractional integrals and Itô integrals, and derive the fractional Itô formula within our framework. We also examine differential equations driven by rough paths, with a focus on the fractional Black-Scholes model, establishing the chain rule for this model and demonstrating the zero expectation property of our Itô integration. In Section \ref{sec-general-integration}, we broaden our discussion to more general Itô integrals beyond the simple one-form associated with fBM, particularly for fractional Ornstein-Uhlenbeck processes. We also consider the integration for solutions to rough differential equations, discussing how renormalization is contingent upon the nature of the integrand. Section \ref{Finance-sec} applies our developed theories to financial markets, illustrating the absence of arbitrage in the Itô fractional Black-Scholes market under specific trading strategy restrictions. Finally, Section \ref{sec-fOU-estimate} focuses on the application of our integration theory to the estimation problems in fOU processes.

As a note, the development of a general pathwise Itô integration theory tailored to fractional Brownian motion necessitates a specialized approach to the Itô rough path lift of fBM, specifically its Lévy area, which uniquely depends on the integrand. This requirement starkly contrasts with the classical Itô's theory for Brownian motion. In this paper, we will demonstrate that Itô integrals with respect to fBM can be constructed using a universal rough path lift of fBM, applicable across all one-forms. However, when the integrand itself is a solution to a rough differential equation driven by fBM, the specific Itô rough path lift must also adapt to accommodate the solution's dynamics. This approach bears resemblance to the methodology in regularity structures, where the models are intricately dependent on the underlying dynamics.

\section{Preliminaries on rough paths and fractional Brownian motion}\label{sec-prelimi}
\subsection{Rough paths}
\label{RPT-sec} 
This section aims to establish several notations that will be used throughout the paper by recalling some basic notions concerning rough paths. Our exposition closely follows the rough path literature (e.g., \cite{Allan21,CQ02,FH14,FV10,LQ02}). We will focus on the fundamental framework required to ensure that fBM with Hurst parameter $H\in(\frac{1}{3},\frac{1}{2}]$, which we introduce in the next subsection, has a natural rough path lift.

For $N\in\mathbb{N}$, $T^{(N)}(\mathbb{R}^{d})$ denotes the truncated
tensor algebra defined by
\[
T^{(N)}(\mathbb{R}^{d}):=\oplus_{n=0}^{N}(\mathbb{R}^{d})^{\otimes n},
\]
with the convention that $(\mathbb{R}^{d})^{\otimes0}=\mathbb{R}.$
The space $T^{(N)}(\mathbb{R}^{d})$ is equipped with a vector space
structure and a multiplication $\otimes$ defined by
\[
(\bfX\otimes \bfY)^{k}=\sum_{i=0}^{k}\bfX^{k-i}\bfY^{i},\ \ k=0,1,\cdots,N,
\]
where $\bfX=(1,\bfX^{1},\cdots,\bfX^{N})$, $\bfY=(1,\bfY^{1},\cdots,\bfY^{N})\in T^{(N)}(\mathbb{R}^{d}).$

We will consider continuous $\mathbb{R}^{d}$-valued paths $X$ on
$[0,T]$ with bounded variations, and their canonical lifts $\bfX_{s,t}=(1,\bfX_{s,t}^{1},\cdots,\bfX_{s,t}^{N})$
in the space $T^{(N)}(\mathbb{R}^{d})$, where
\[
\bfX_{s,t}^{1}=X_{t}-X_{s},
\]
\[
\bfX_{s,t}^{2}=\int_{s<t_{1}<t_{2}<t}dX_{t_{1}}\otimes dX_{t_{2}},
\]
\[
\vdots
\]
and
\[
\bfX_{s,t}^{N}=\int_{s<t_{1}<\cdots<t_{N}<t}dX_{t_{1}}\otimes\cdots\otimes dX_{t_{N}}.
\]
The lifted path satisfies \textit{``Chen's identity"} as following:
\begin{equation}
\bfX_{s,t}=\bfX_{s,u}\otimes \bfX_{u,t},\ \forall(s,u),\ (u,t)\in\Delta,\label{Chen}
\end{equation}
where $\Delta$ denotes the simplex $\{(s,t):0\leq s<t\leq T\}$.

By definition, a continuous map $\bfX$ from the simplex $\Delta$ into
a truncated tensor algebra $T^{(N)}(\mathbb{R}^{d})$ is called a
\textit{rough path} (of roughness $p\geq1$, where $N=[p]$), if it
satisfies (\ref{Chen}) and has finite $p$-variations, that is,
\[
\sum_{i=1}^{N}\sup_{D}\sum_{\ell}|\bfX_{t_{\ell-1},t_{\ell}}^{i}|^{p/i}<\infty,
\]
where the sup runs over all finite partitions $D=\{0=t_{0}<t_{1}<\cdots<t_{n}=T\}$.
The $p$-variation distance is defined to be
\[
d_{p}(\bfX,\bfY)=\sum_{i=1}^{N}\left(\sup_{D}\sum_{\ell}|\bfX_{t_{\ell-1},t_{\ell}}^{i}-\bfY_{t_{\ell-1},t_{\ell}}^{i}|^{p/i}\right)^{i/p}.
\]

Equivalently, $\bfX:\Delta\to T^{(N)}(\mathbb{R}^{d})$ has finite $p$-variations
if
\[
|\bfX_{s,t}^{i}|\leq\omega(s,t)^{i/p},\ \forall i=1,\cdots,N,\ \forall(s,t)\in\Delta
\]
for some function $\omega$, where $\omega$ is a non-negative, continuous,
super-additive function on $\Delta$ and satisfies $\omega(t,t)=0$.
Such function $\omega$ is called a \textit{control of the rough path
$\bfX$.}

The space of all $p$-rough paths is denoted by $\Omega_{p}(\mathbb{R}^{d})$.
A rough path $\bfX$ is called a geometric rough path if there is a sequence
of $\bfX(n)$, where $\bfX(n)$ are the canonical lifts of their first level
$\bfX(n)^{1}$ which are continuous with finite variations, such that
$\bfX$ is the limit of $\bfX(n)$ under $p$-variation distance $d_{p}$.
The space of geometric rough paths is denoted by $G\Omega_{p}(\mathbb{R}^{d})$.

As our interest lies in fBMs with Hurst parameter $H\in(\frac{1}{3},\frac{1}{2}]$,
which will be introduced later, we consider only rough paths valued
in $T^{(2)}(\mathbb{R}^{d})$. Thus, in what follows, we will assume
that $3>p\geq2$ so that $[p]=2$. A rough path of roughness $p$
can be written as $\bfX_{s,t}=(1,\bfX_{s,t}^{1},\bfX_{s,t}^{2})$ for $s<t$,
and the algebraic relation (Chen's identity) now becomes
\begin{equation}
\bfX_{s,t}^{1}=X_{t}-X_{s},\label{Chen-r1}
\end{equation}
and
\begin{equation}
\bfX_{s,t}^{2}-\bfX_{s,u}^{2}-\bfX_{u,t}^{2}=\bfX_{s,u}^{1}\otimes \bfX_{u,t}^{1},\label{Chen-r2}
\end{equation}
for all $(s,u),\ (u,t)\in\Delta$, where $\bfX^{2}$ should be (if it makes
sense) considered as an iterated integral
\begin{equation}
\int_{s}^{t}\bfX_{s,u}^{1}d\bfX_{u}^{1}:=\bfX_{s,t}^{2}
\end{equation}
which is of course not defined a priori.

The most convenient tool to construct rough paths is through \textit{almost
rough path}s. A function $\bfY=(1,\bfY^{1},\bfY^{2})$ from $\Delta$ to $T^{(2)}(\mathbb{R}^{d})$
is called an \textit{almost rough path} if it has finite $p$-variation,
and for some control $\omega$ and constant $\theta>1$,
\begin{equation}
|(\bfY_{s,t}\otimes \bfY_{t,u})^{i}-\bfY_{s,u}^{i}|\leq\omega(s,u)^{\theta},\ i=1,2,
\end{equation}
for all $(s,t),(t,u)\in\Delta$. According to Theorem 3.2.1 in Lyons and Qian \cite{LQ02},
given an almost rough path $\bfY=(1,\bfY^{1},\bfY^{2})$, there exists a unique
rough path $\bfX=(1,\bfX^{1},\bfX^{2})$ such that
\begin{equation}
|\bfX_{s,t}^{i}-\bfY_{s,t}^{i}|\leq\omega(s,t)^{\theta},\ i=1,2,\ \theta>1
\end{equation}
for some control $\omega$, and all $(s,t)\in\Delta$. 

As a clarification, unless otherwise indicated, we may utilize the symbol $\bfX=(1,X,\bbX)$ as an alternative representation of $\bfX=(1,\bfX^{1},\bfX^{2})$, provided this does not introduce any ambiguity.

\subsection{Fractional Brownian motion}\label{fBM-sec}

Fractional Brownian motion (fBM) $B^H(t)$ is a continuous-time Gaussian process, where $t\geq 0$. It has a mean of zero for all $t\geq 0$, and the covariance function is defined by the following equation:
\begin{equation}
E[B^{H}(t)B^{H}(s)]={\frac{1}{2}}(|t|^{2H}+|s|^{2H}-|t-s|^{2H}),
\end{equation}
where $H\in(0,1)$ is the Hurst parameter.

When \( H>1/2 \), the increments of fractional Brownian motion (fBM) exhibit positive correlation, indicating persistence: a tendency for future increments to continue in the same direction as past increments. Conversely, for \( H<1/2 \), the increments are negatively correlated, demonstrating counter-persistence: a tendency for future increments to reverse the direction of past increments. Thus, when \( H<1/2 \), fBM is likely to show a decrease following an increase and vice versa. This dual characteristic makes fBM a versatile model for phenomena exhibiting both short-range and long-range dependencies across various disciplines, including physics, biology, hydrology, network research, and financial mathematics.

The fBM with Hurst parameter $H$ has an integral representation in
terms of Brownian motion
\begin{equation}
B^{H}(t)=\int_{0}^{t}K_{H}(t,s)dW(s),
\end{equation}
where $W(t)$ is a standard Brownian motion and
\[
K_{H}(t,s)=C_{H}\left[\frac{2}{2H-1}\left(\frac{t(t-s)}{s}\right)^{H-\frac{1}{2}}
-\int_{s}^{t}\left(\frac{u(u-s)}{s}\right)^{H-\frac{1}{2}}\frac{du}{u}\right]1_{(0,t)}(s)
\]
which is a singular kernel, and $C_{H}$ is a normalised constant.

Integration theories for fractional Brownian motion with a Hurst parameter \( H>\frac{1}{2} \) can be developed using Young's integration theory or functional integration approaches, as outlined in works such as \cite{CKR93,Zah98}. On the other hand, stochastic calculus for fBM with \( H<\frac{1}{2} \) is more suitably addressed within the framework of rough path analysis. In their seminal work, Coutin and Qian \cite{CQ02} describe the construction of a canonical level-3 rough path \( B^{H} \) for fBM with \( H>1/4 \). This construction includes defining iterated integrals of multi-dimensional fBM up to level-3, enabling the application of rough path integration theory to effectively handle fBM.

The second and third level processes are defined in terms of iterated
Riemann-Stieltjes integrals along the dyadic piece-wise linear approximations,
and geometric rough paths of fBM are the limits in $p$-variation
distance. Here as our main concern is the fBM with Hurst parameter
$H\in(\frac{1}{3},\frac{1}{2}]$, so we only need the level-2 results
about fBM. The method of \cite{CQ02} to construct the fractional
Brownian rough path $\bfB^{H}=(1,B^{H},\bbB^{H})$ implies that
\begin{equation}
\bbB^{H}_{s,t}:=\lim_{m\to0}\int_{s}^{t}B_{s,u}^{H,(m)}\otimes dB_{u}^{H,(m)},\ \ a.s.,\label{fBm-rough}
\end{equation}
exists in $p$-variation distance as long as $pH>1$, where the dyadic
piece-wise linear approximations $B_{t}^{H,(m)}$ on interval $[s,t]$
is defined by
\[
B_{r}^{H,(m)}:=B_{t_{\ell-1}^{m}}^{H}+2^{m}\frac{r-t_{\ell-1}^{m}}{t_{\ell}^{m}-t_{\ell-1}^{m}}\Delta_{\ell}^{m}B^{H},
\]
with $\Delta_{\ell}^{m}B^{H}=B_{t_{\ell}^{m}}^{H}-B_{t_{\ell-1}^{m}}^{H}$,
$t_{\ell}^{m}:=s+\frac{\ell}{2^{m}}(t-s)$ for $\ell=1,2,\cdots,2^{m}$.

\begin{proposition}(\cite{CQ02}, Theorem 2; \cite{FH14}, Theorem 10.4) Let $B^{H}=(B^{H,1},\cdots,B^{H,d})$ be a $d$-dimensional fBM with
the Hurst parameter $H\in(\frac{1}{3},\frac{1}{2}]$. Then $B^{H}$,
restricted to an finite interval $[0,T]$, lifts via (\ref{fBm-rough})
to a geometric rough path $\bfB^{H}=(1,B^{H},\bbB^{H})\in G\Omega_{p}([0,T],\mathbb{R}^{d})$,
for all $p\in(\frac{1}{H},3)$.

\end{proposition}

The proof details for this proposition are documented in \cite{CQ02,FH14}. Here, the random rough path $\bfB^{H}=(1,B^{H},\bbB^{H})$ is referred to as the canonical lift, which parallels the Stratonovich lift of fBM. In this paper, however, we introduce a novel natural lift of fBM in the Itô sense, prompting us to designate the Stratonovich fractional Brownian rough path $\bfB^{H}=(1,B^{H},\bbB^{H})$ as $\bfB^{H,\text{Str}}=(1,B^{H},\bbB^{H,\text{Str}})$. For brevity, we simply use $B$ to denote $B^{H}$ when there is no risk of confusion.

We define the Itô rough path associated with fBM $B^{H}$ as follows:
\[
\bfB^{H,\text{Itô}}_{s,t}=(1,B^{H}_{s,t},\bbB^{H,\text{Itô}}_{s,t})=\left(1,B_{t}^H-B_{s}^H,\bbB^{H,\text{Str}}_{s,t}-\varphi_{s,t}\right),
\]
where the Hurst parameter \( H \) is within \( (\frac{1}{3},\frac{1}{2}] \). It is established that $\bfB^{H,\text{Itô}}=(1,B^{H},\bbB^{H,\text{Itô}})$ constitutes a random rough path, albeit a non-geometric one. We term this non-geometric rough path the Itô fractional Brownian rough path, as the rough path and its integration theory extend the traditional frameworks of standard Brownian motion and Itô stochastic integration. To simplify notation, we denote it by $\bfB_{s,t}=(1,B_{s,t},\bbB_{s,t})$ where appropriate to avoid confusion.

\section{It\^o Integration against fractional Brownian motion}

\label{RI-sec}

\subsection{It\^o integrals of one forms against fBM}

The objective of this section is to define Itô integrals of one-forms with respect to fractional Brownian motion, exemplified by expressions such as \(\int_{s}^{t}F(B) \, dB\), where \(B=(B^{(1)}, \ldots, B^{(d)})\) represents a \(d\)-dimensional fBM characterized by a Hurst parameter \(H\) in the range \((\frac{1}{3}, \frac{1}{2}]\).

In order to define $\int_{s}^{t}F(B)dB$, where $F:\mathbb{R}^{d}\to L(\mathbb{R}^{d},\mathbb{R}^{e})$
satisfies some smoothness conditions, according to the rough path theory
of Lyons, we should take $B$ as a rough path. Actually the symbol
$\int_{s}^{t}F(B)dB$, to some extent, is a misleading statement.
Recall that the rough integral $\int F(X)d\bfX$ against a rough path
$\bfX=(1,X,\bbX)\in\Omega_{p}(\mathbb{R}^{d})$ with $2\leq p<3$
is defined to be again a rough path. The rough integral is defined
uniquely by its associated almost rough path $\widehat{\bfY}=(1,\widehat{Y},\widehat{\bbY})$ as following (see e.g. Definition 5.2.1 in Lyons and Qian~\cite{LQ02}):
\begin{align}
\widehat{Y}_{s,t} & =F(X_{s})X_{s,t}+DF(X_{s})\bbX_{s,t},\\
\widehat{\bbY}_{s,t} & =F(X_{s})\otimes F(X_{s})\bbX_{s,t},
\end{align}
and the integral $\int F(X)d\bfX$ is defined to be the rough path $\bfY=(1,Y,\bbY)$ uniquely associated with the almost rough path $\widehat{\bfY}=(1,\widehat{Y},\widehat{\bbY})$ (see e.g. Theorem 5.2.1 in \cite{LQ02}). The integral can be written
in compensated Riemann sum form, that is,
\begin{equation}
\begin{split}Y_{s,t} =\int_{s}^{t}F(X)\rd\bfX  :=\lim_{|D|\to0}\sum_{\ell}\left(F(X_{t_{\ell-1}})X_{t_{\ell-1},t_{\ell}}+DF(X_{t_{\ell-1}})\bbX_{t_{\ell-1},t_{\ell}}\right)
\end{split}
\end{equation}
and the second level
\begin{equation}
\begin{split}\bbY_{s,t} =\int_{s}^{t}F(X)\rrd\bfX:=\lim_{|D|\to0}\sum_{\ell}\left(Y_{s,t_{\ell-1}}\otimes Y_{t_{\ell-1},t_{\ell}}+F(X_{t_{\ell-1}})\otimes F(X_{t_{\ell-1}})\bbX_{t_{\ell-1},t_{\ell}}\right).
\end{split}
\end{equation}
These limits exist in the $p$-variation distance.

Using the general definition provided earlier, we calculate the integral against the Stratonovich fractional Brownian rough path $\bfB^{H,\text{Str}}$, denoted as $\bfX^{\text{Str}} := \int F(B^H) \circ d\bfB^{H,\text{Str}}$. For clarity and simplicity, we use $\bfX^{\text{Str}} = \int F(B) \circ d\bfB$ to represent this Stratonovich integral, which is defined against the rough path $\bfB^{H,\text{Str}}$.

Similarly, we define the Itô integral $\int F(B^H) d\bfB^{H,\text{Itô}}$ against the rough path $\bfB^{H,\text{Itô}}$, where $\varphi(t) = \frac{1}{2} t^{2H}$. This integral is denoted as $\bfX^{\text{Itô}} := \int F(B) d\bfB$, simplifying the notation for practical purposes.

\subsubsection{Relation between Stratonovich and It\^o rough integrals (I)}

Now we establish a relation between Stratonovich and It\^o integrals.

\begin{theorem}\label{SI-ingtral-thm-notime} The relation between Stratonovich and It\^o integral is
given as the following.\\
 (i) For the first level,
\begin{equation}
X^{\text{Str}}_{s,t}-X^{\text{It\^o}}_{s,t}=\frac{1}{2}\int_{s}^{t}DF(B_{u})du^{2H},
\end{equation}
(ii) For the second level,
\begin{equation}
\begin{split}\bbX^{\text{Str}}_{s,t}-\bbX^{\text{It\^o}}_{s,t} & =\frac{1}{2}\int_{s}^{t}F(B_{u})\otimes F(B_{u})du^{2H}
 +\frac{1}{2}\int_{s}^{t}\left(\int_{s}^{u}DF(B_{r})dr^{2H}\right)\otimes dX^{\text{Str}}_{0,u}\\
 & +\frac{1}{2}\int_{s}^{t}X^{\text{Str}}_{s,u}\otimes DF(B_{u})du^{2H} -\frac{1}{4}\int_{s}^{t}\left(\int_{s}^{u}DF(B_{r})dr^{2H}\right)DF(B_{u})du^{2H},
\end{split}
\label{SI-2}
\end{equation}
where the last four integrals are Young integrals.
\end{theorem}

See Appendix~\ref{appendix_SI_correction} for proof. The Stratonovich-It\^o correction follows from more general results on the relation between geometric and non-geometric lifts, see e.g. Theorem 5.4.2 in Lyons and Qian~\cite{LQ02} with $\phi_t = -t^{2H}$ there.

\subsubsection{Relation between Stratonovich and It\^o rough integrals (II)}

\noindent  Let us introduce the space-time Stratonovich/It\^o path
$\widetilde{B}=(B,t)$, where the first level is given by
\[
\widetilde{B}_{s,t}=(B_{s,t},t-s),
\]
and the second level is given by
\[
\widetilde{\bbB}_{s,t}=\left(\bbB_{s,t},\int_{s}^{t}B_{s,u}du,\int_{s}^{t}(u-s)dB_{u},\frac{1}{2}(t-s)^{2}\right),
\]
where the cross integrals are Young integrals, and $\bbB$ is the second level Stratonovich or It\^o lift of fBM. Naturally
\begin{equation}
\int F(B,t)d\bfB:=\int f(\widetilde{B})d\widetilde{\bfB},
\end{equation}
with $f(x,t)(\xi,\tau)=F(x,t)\xi$ and the right hand side is well
defined as an integral for rough paths. We use the symbol $\int F(B,t)d\bfB=:\bfX^{\text{It\^o}}$
as It\^o integral and the symbol $\int F(B,t)\circ d\bfB=:\bfX^{\text{Str}}$ as Stratonovich
integral. We can also establish the relationship between the Stratonovich and It\^o integrals.

\begin{theorem}\label{SI-ingtral-thm} The relation between Stratonovich
and It\^o integral is given as the following.\\
 (i) For the first level,
\begin{equation}
X^{\text{Str}}_{s,t}-X^{\text{It\^o}}_{s,t}=\frac{1}{2}\int_{s}^{t}D_{x}F(B_{u},u)du^{2H},\label{SI-11}
\end{equation}
(ii) For the second level,
\begin{equation}
\begin{split}\bbX^{\text{Str}}_{s,t}-\bbX^{\text{It\^o}}_{s,t} & =\frac{1}{2}\int_{s}^{t}F(B_{u},u)\otimes F(B_{u},u)du^{2H} +\frac{1}{2}\int_{s}^{t}\left(\int_{s}^{u}D_{x}F(B_{r})dr^{2H}\right)\otimes dX^{\text{Str}}_{0,u}\\
 & +\frac{1}{2}\int_{s}^{t}X^{\text{Str}}_{s,u}\otimes D_{x}F(B_{u})du^{2H} -\frac{1}{4}\int_{s}^{t}\left(\int_{s}^{u}D_{x}F(B_{r},r)dr^{2H}\right)D_{x}F(B_{u},u)du^{2H},
\end{split}
\label{SI-22}
\end{equation}
where the last four integrals are Young integrals. 
\end{theorem}

The proof can be found in Appendix~\ref{appendix_SI_correction}. The inhomogeneous case follows with the space-time lift.  See also, for example, Theorem 5.4.2 in Lyon and Qian~\cite{LQ02}.

\subsection{The It\^o formula}

\subsubsection{Homogeneous It\^o type formula}

Let $\bfX=(1,X,\bbX)\in\Omega_{p}(\mathbb{R}^{d})$, $2\leq p<3$
be a $p$-rough path, and $F:\mathbb{R}^{d}\to L(\mathbb{R}^{d},\mathbb{R}^{e})$
be a $\mathrm{Lip}(\gamma)$ function for some $\gamma>p$ (see Definition 5.1.1 in \cite{LQ02}). Since often composes with rough paths, we want to make the function $F(X)$ into a rough path $F_{\mathfrak{R}}(\bfX)=(1,F_{\mathfrak{R}}(\bfX)^{1},F_{\mathfrak{R}}(\bfX)^{2})$.
In terms of rough path integrals, we use the formula
\begin{equation}
F_{\mathfrak{R}}(\bfX)=\int DF(X)d\bfX
\end{equation}
as a definition of image $F_{\mathfrak{R}}$ of function $F$. Actually
$F_{\mathfrak{R}}$ sends a rough path to another rough path, so that
$F_{\mathfrak{R}}:\Omega_{p}(\mathbb{R}^{d})\to\Omega_{p}(\mathbb{R}^{e})$,
$2\leq p<3$. By the definition of rough path integrals,
\begin{equation}
\begin{split}F_{\mathfrak{R}}(\bfX)_{s,t}^{1} & =\int_{s}^{t}DF(X)\rd\bfX  =\lim_{|D|\to0}\sum_{\ell}\left(DF(X_{t_{\ell-1}})X_{t_{\ell-1},t_{\ell}}+D^{2}F(X_{t_{\ell-1}})\bbX_{t_{\ell-1},t_{\ell}}\right),
\end{split}
\end{equation}
\begin{equation}
\begin{split}F_{\mathfrak{R}}(\bfX)_{s,t}^{2} & =\int_{s}^{t}DF(X)\rrd\bfX\\
 & =\lim_{|D|\to0}\sum_{\ell}\left(F_{\mathfrak{R}}(\bfX)_{s,t_{\ell-1}}^{1}\otimes F_{\mathfrak{R}}(\bfX)_{t_{\ell-1},t_{\ell}}^{1} +DF(X_{t_{\ell-1}})\otimes DF(X_{t_{\ell-1}})\bbX_{t_{\ell-1},t_{\ell}}\right).
\end{split}
\end{equation}

Let $\bfX(\phi)_{s,t}=(1,\bfX(\phi)_{s,t}^{1},\bfX(\phi)_{s,t}^{2}):=(1,X_{s,t},\bbX_{s,t}-\phi_{s,t})$
be a perturbation of the rough path $\bfX=(1,X,\bbX)$, where $\phi_{s,t}=\phi_{t}-\phi_{s}$
(additive) is finite $q$-variation with $q\leq\frac{p}{2}$. Assume
that $\bfX$ is a geometric rough path, then $\bfX(\phi)$ is no longer
a geometric rough path in general. Define the composition $F_{\mathfrak{R}}(\bfX(\phi))=(1,F_{\mathfrak{R}}(\bfX(\phi))^{1},F_{\mathfrak{R}}(\bfX(\phi))^{2})$
as
\begin{equation}
F_{\mathfrak{R}}(\bfX(\phi))=\int DF(X(\phi))d\bfX(\phi).
\end{equation}
We have the following It\^o type formula.

\begin{theorem}\label{Ito_FX}(It\^o type formula) Assume that $\bfX\in G\Omega_{p}(\mathbb{R}^{d})$
with $2\leq p<3$, $\bfX(\phi)$ is a perturbation of the rough path
$\bfX$ as above, and $F:\mathbb{R}^{d}\to L(\mathbb{R}^{d},\mathbb{R}^{e})$
is a $\mathrm{Lip}(\gamma)$ function for some $\gamma>p$, then\\
 (i) $F_{\mathfrak{R}}(\bfX)_{s,t}^{1}=F(X_{t})-F(X_{s})$.\\
 (ii) For the first level,
\[
F(X_{t})-F(X_{s})=\int_{s}^{t}DF(X(\phi))\rd\bfX(\phi)+\int_{s}^{t}D^{2}F(X_{u})d\phi_{u}.
\]
(iii) For the second level,
\[
\begin{split}F_{\mathfrak{R}}(\bfX)_{s,t}^{2} & =\int_{s}^{t}DF(X(\phi))\rrd\bfX(\phi)+\int_{s}^{t}DF(X_{u})\otimes DF(X_{u})d\phi_{u} +\int_{s}^{t}\left(\int_{s}^{u}D^{2}F(X_{r})d\phi_{r}\right)\otimes dF(X_{u})\\
 & \ \ \ +\int_{s}^{t}(F(X_{u})-F(X_{s}))\otimes D^{2}F(X_{u})d\phi_{u} +\int_{s}^{t}\left(\int_{s}^{u}D^{2}F(X_{r})d\phi_{r}\right)D^{2}F(X_{u})d\phi_{u},
\end{split}
\]
where the last four integrals are Young integrals, and $F_{\mathfrak{R}}(\bfX)_{s,t}^{2}$
can be viewed as a kind of geometric increments of the second level
process.

\end{theorem}

\begin{proof} (i) The equality holds for any continuous path $X$
with finite variation and its canonical lift as rough paths. Then
by definition of geometric rough paths, it can be approximated by
a sequence of path with finite variation in $p$-variation. By continuity
of both sides, we know the equality still holds.
 (ii) and (iii) can be proved by the same arguments as in Theorem
\ref{SI-ingtral-thm}. A proof can also be found in Friz and Hairer~\cite{FH14} as the proof of Equation (5.5) there.
\end{proof}

\begin{remark} In integration theory for rough paths, if $\bfX\in\Omega_{p}$,
$2\leq p<3$, one cannot just write a symbol
$dF(X_{t})$ as the ordinary case. There is no meaning for this symbol
unless under the sense of Young integrals, if it is well defined.
We can see an example below which says that the $dF(X_{t})$ is undefined
in general. Actually if we want to make sense the differential symbol,
we should lift $F(X_{t})$ to a rough path $F_{\mathfrak{R}}(\bfX)$
as above and then understand the differential symbol as
\[
dF_{\mathfrak{R}}(\bfX)=d(F_{\mathfrak{R}}(\bfX)^{1},F_{\mathfrak{R}}(\bfX)^{2}).
\]
\end{remark}

To clarify the above remark, we first give a lemma below.

\begin{lemma}\label{dY} Let $\bfY=(1,Y,\bbY)$ be a rough path.
Then the integral
\begin{equation}
\int_{s}^{t}d\bfY=(1,Y_{s,t},\bbY_{s,t})
\end{equation}
as expected.

\end{lemma}

\noindent\textbf{Example 1.} Take $\bfY$ as $F_{\mathfrak{R}}(\bfX)$, $F_{\mathfrak{R}}(\bfX(\phi))$.
Then by lemma \ref{dY}, we have
\[
\int_{s}^{t}dF_{\mathfrak{R}}(\bfX)=(1,F_{\mathfrak{R}}(\bfX)_{s,t}^{1},F_{\mathfrak{R}}(\bfX)_{s,t}^{2}),
\]
\[
\int_{s}^{t}dF_{\mathfrak{R}}(\bfX(\phi))=(1,F_{\mathfrak{R}}(\bfX(\phi))_{s,t}^{1},F_{\mathfrak{R}}(\bfX(\phi))_{s,t}^{2}).
\]
Actually $F(X(\phi)_{t})=F(X_{t})$, but $dF_{\mathfrak{R}}(\bfX)\neq dF_{\mathfrak{R}}(\bfX(\phi))$,
even for the first level as we can see that $F_{\mathfrak{R}}(\bfX)_{s,t}^{1}\neq F_{\mathfrak{R}}(\bfX(\phi))_{s,t}^{1}$
by It\^o formula above. Therefore the symbol $dF(X_{t})$ or $dF(X(\phi)_{t})$
for rough paths can lead to confusion.

\begin{remark} If $X$ is a geometric rough path, by Theorem \ref{Ito_FX},
we have
\[
F_{\mathfrak{R}}(\bfX)_{s,t}^{1}=F(X_{t}^{1})-F(X_{s}^{1}),
\]
i.e. $F_{\mathfrak{R}}(\bfX)_{s,t}^{1}=F(X^{1})_{s,t}$. However, for
the non-geometric rough path we do not have the equality alike. In
fact, in general,
\[
F_{\mathfrak{R}}(\bfX(\phi))_{s,t}^{1}\neq F(X(\phi)_{t}^{1})-F(X(\phi)_{s}^{1})(=F(X_{t}^{1})-F(X_{s}^{1})).
\]
 \end{remark}

\noindent We next want to establish an It\^o type formula for integrals
against It\^o fractional Brownian rough path. As a corollary, we have

\noindent \begin{corollary}(It\^o formula for fBM) Let $\bfB^{\text{Str}}$ be fractional
Brownian rough path with Hurst parameter $\frac{1}{3}<H\leq\frac{1}{2}$
enhanced under Stratonovich sense, $\bfB^{\text{It\^o}}$ be the It\^o fractional Brownian
rough path, and $F:\mathbb{R}^{d}\to L(\mathbb{R}^{d},\mathbb{R}^{e})$
be a $\mathrm{Lip}(\gamma)$ function for some $H\gamma>1$. Then\\
 (i) $F(B_{t})-F(B_{s})=F_{\mathfrak{R}}(\bfB^{\text{Str}})_{s,t}^{1}=\int_{s}^{t}DF(B)\circ \rd\bfB^{\text{Str}}$.\\
 (ii) For the first level,
\[
F(B_{t})-F(B_{s})=\int_{s}^{t}DF(B)\rd\bfB^{\text{It\^o}}+\frac{1}{2}\int_{s}^{t}D^{2}F(B_{u})du^{2H}.
\]
(iii) For the second level,
\[
\begin{split}F_{\mathfrak{R}}(\bfB^{\text{Str}})_{s,t}^{2} & =\int_{s}^{t}DF(B)\rrd\bfB^{\text{It\^o}}+\frac{1}{2}\int_{s}^{t}DF(B_{u})\otimes DF(B_{u})du^{2H} +\frac{1}{2}\int_{s}^{t}\left(\int_{s}^{u}D^{2}F(B_{r})dr^{2H}\right)\otimes dF(B_{u})\\
 & \ \ \ +\frac{1}{2}\int_{s}^{t}(F(B_{u})-F(B_{s}))\otimes D^{2}F(B_{u})du^{2H} -\frac{1}{4}\int_{s}^{t}\left(\int_{s}^{u}D^{2}F(B_{r})dr^{2H}\right)\otimes D^{2}F(B_{u})du^{2H},
\end{split}
\]
where the last four integrals are Young integrals. \end{corollary}

\subsubsection{Inhomogeneous It\^o formula}

In the following, we want to make sense of $F_{\mathfrak{R}}((\bfX,t))$
when the inhomogeneous function $F(x,t)$ is applied to a rough path
$\bfX=(1,X,\bbX)\in\Omega_{p}$ and establish It\^o formula for it.
First, recall the space-time rough path $\widetilde{\bfX}=(\bfX,t)$, where
the first level is given by
\[
\widetilde{X}_{s,t}=(X_{s,t},t-s),
\]
and the second level is given by
\[
\widetilde{\bbX}_{s,t}=\left(\bbX_{s,t},\int_{s}^{t}X_{s,u}du,\int_{s}^{t}(u-s)dX_{u},\frac{1}{2}(t-s)^{2}\right),
\]
where the cross integrals are Young integrals. Define the rough path
$F_{\mathfrak{R}}((\bfX,t))$ by
\begin{equation}
F_{\mathfrak{R}}((\bfX,t)):=F_{\mathfrak{R}}(\widetilde{\bfX})=\int DF(\widetilde{X})d\widetilde{\bfX},
\end{equation}
where $DF(x,t)(\xi,\tau)=D_{x}F(x,t)\xi+D_{t}F(x,t)\tau$.

Note that if $\bfX(\phi)$ is a perturbation of the rough path $\bfX$,
and $\widetilde{\bfX}(\phi)$ is its associated space-time rough path,
then
\[
\widetilde{X}(\phi)_{s,t}=\widetilde{X}_{s,t},
\]
\[
\widetilde{\bbX}(\phi)_{s,t}=\left(\bbX_{s,t}-\phi_{s,t},\int_{s}^{t}X_{s,u}du,\int_{s}^{t}(u-s)dX_{u},\frac{1}{2}(t-s)^{2}\right).
\]
Note that only the first $d\times d$ dimensional components of the
second level of $\widetilde{\bfX}$ are changed.

\begin{theorem}\label{Ito_Ft}(It\^o formula) Assume $\bfX\in G\Omega_{p}(\mathbb{R}^{d})$
with $2\leq p<3$, $\bfX(\phi)$ is a perturbation of the rough path
$\bfX$, and $\widetilde{\bfX}$, $\widetilde{\bfX}(\phi)$ are their associated
space-time rough path respectively, $F:\mathbb{R}^{d+1}\to L(\mathbb{R}^{d+1},\mathbb{R}^{e})$
be a $\mathrm{Lip}(\gamma)$ function for some $\gamma>p$.\\
 (i) We have the basic calculus formula:
\begin{equation}
F(X_{t},t)-F(X_{s},s)=\int_{s}^{t}DF(\widetilde{X})\rd\widetilde{\bfX}.\label{Ito_Stratcase}
\end{equation}
(ii) For the first level,
\begin{equation}
F(X_{t},t)-F(X_{s},s)=\int_{s}^{t}DF(\widetilde{X}(\phi))\rd\widetilde{\bfX}(\phi)+\int_{s}^{t}D_{x}^{2}F(X_{u},u)d\phi_{u}.\label{Ito_Ft}
\end{equation}
(iii) For the second level,
\[
\begin{split}
&\quad F_{\mathfrak{R}}((\bfX,t))_{s,t}^{2} \\
&=\int_{s}^{t}DF(\widetilde{X}(\phi))\rrd\widetilde{\bfX}(\phi)+\int_{s}^{t}D_{x}F(X_{u},u)\otimes D_{x}F(X_{u},u)d\phi_{u} +\int_{s}^{t}\left(\int_{s}^{u}D_{x}^{2}F(X_{r},r)d\phi_{r}\right)\otimes dF(X_{u},u)\\
 & \ \ \ +\int_{s}^{t}(F(X_{u},u)-F(X_{s},s))\otimes D_{x}^{2}F(X_{u},u)d\phi_{u} +\int_{s}^{t}\left(\int_{s}^{u}D_{x}^{2}F(X_{r},r)d\phi_{r}\right)\otimes D_{x}^{2}F(X_{u},u)d\phi_{u},
\end{split}
\]
where the last four integrals are Young integrals.

\end{theorem}

\begin{proof} (i) The proof is same as (i) of Theorem \ref{Ito_FX},
first for $p=1$ it holds, then the result for any geometric rough
path follows from the continuity.
 For (ii) and (iii), note that
\[
\left|\int_{s}^{t}X_{s,u}du\right|=o(|t-s|)
\]
and
\[
\left|\int_{s}^{t}(u-s)\rd\bfX_{u}\right|=o(|t-s|),
\]
the rest of proof is same as Theorem \ref{SI-ingtral-thm}.
\end{proof}

Note that if $X_{t}$ is a continuous path with
finite variation, then Equation (\ref{Ito_Stratcase}) reads as
\begin{equation}
F(X_{t},t)-F(X_{s},s)=\int_{s}^{t}D_{x}F(X_{u},u)dX_{u}+\int_{s}^{t}D_{u}F(X_{u},u)du,
\end{equation}
and Equation (\ref{Ito_Ft}) becomes
\begin{equation}
\begin{split}F(X_{t},t)-F(X_{s},s) & =\int_{s}^{t}D_{x}F(X_{u},u)dX_{u}+\int_{s}^{t}D_{u}F(X_{u},u)du +\int_{s}^{t}D_{x}^{2}F(X_{u},u)d\phi_{u}.
\end{split}
\end{equation}
These equations are just like It\^o formulas in terms of the Stratonovich
and It\^o integrals in stochastic calculus.

Now set $\widetilde{\bfB}^{\text{Str}}$, $\widetilde{\bfB}^{\text{It\^o}}$ are the associated
space-time rough paths of Stratonovich fractional Brownian rough path
$\bfB^{\text{Str}}$ and It\^o rough path $\bfB^{\text{It\^o}}$, respectively.

\begin{corollary}(It\^o formula for fractional Brownian rough path) Let $\bfB^{\text{Str}}$
be fractional Brownian rough path with Hurst parameter $\frac{1}{3}<H\leq\frac{1}{2}$
enhanced under Stratonovich sense, $\bfB^{\text{It\^o}}$ be the It\^o fractional
Brownian rough path, and $F:\mathbb{R}^{d+1}\to L(\mathbb{R}^{d+1},\mathbb{R}^{e})$
be a $\mathrm{Lip}(\gamma)$ function for some $H\gamma>1$, then\\
 (i) $F(B_{t},t)-F(B_{s},s)=F_{\mathfrak{R}}((\bfB^{\text{Str}},t))_{s,t}^{1}=\int_{s}^{t}DF(\widetilde{B})\circ \rd\widetilde{\bfB}^{\text{Str}}$.\\
 (ii) For the first level,
\begin{equation}
F(B_{t},t)-F(B_{s},s)=\int_{s}^{t}DF(\widetilde{B})\rd\widetilde{\bfB}^{\text{It\^o}}+\frac{1}{2}\int_{s}^{t}D_{x}^{2}F(B_{u},u)du^{2H}.\label{time-dept-Ito-1}
\end{equation}
(iii) For the second level,
\begin{equation}
\begin{split}F_{\mathfrak{R}}((\bfB^{\text{Str}},t))_{s,t}^{2} & =\int_{s}^{t}DF(\widetilde{B})\rrd\widetilde{B}^{\text{It\^o}}+\frac{1}{2}\int_{s}^{t}D_{x}F(B_{u},u)\otimes D_{x}F(B_{u},u)du^{2H}\\
 & \ \ \ +\frac{1}{2}\int_{s}^{t}\left(\int_{s}^{u}D_{x}^{2}F(B_{r},r)dr^{2H}\right)\otimes dF(B_{u},u)\\
 & \ \ \ +\frac{1}{2}\int_{s}^{t}(F(B_{u},u)-F(B_{s},s))\otimes D_{x}^{2}F(B_{u},u)du^{2H}\\
 & \ \ \ -\frac{1}{4}\int_{s}^{t}\left(\int_{s}^{u}D_{x}^{2}F(B_{r},r)dr^{2H}\right)\otimes D_{x}^{2}F(B_{u},u)du^{2H},
\end{split}
\label{time-dept-Ito-2}
\end{equation}
where the last four integrals are Young integrals.

\end{corollary}

\subsection{Differential equations driven by rough paths}

\subsubsection{Basics of differential equations}

In this subsection, we summarize key concepts regarding differential equations driven by rough paths, referencing the framework outlined in \cite{LQ02}. We will revisit the definition of these differential equations to provide an understanding for the discussions that follow.

\begin{definition} Let $f:W\to L(V,W)$ be a vector field on $W$. Let
$y_{0}\in W$ and let $\bfX\in\Omega_{p}(V)$, for $2\leq p<3$. Then
we say that a rough path $\bfY\in\Omega_{p}(W)$ is a solution to the
following initial value problem:
\begin{equation}
d\bfY=f(Y)d\bfX,\ Y_{0}=y_{0}\label{RDEXY}
\end{equation}
if there is a rough path $\bfZ\in\Omega_{p}(V\oplus W)$ such that $\pi_{V}(\bfZ)=\bfX$,
$\pi_{W}(\bfZ)=\bfY$ and
\begin{equation}
\bfZ=\int\widehat{f}(Z)d\bfZ,
\end{equation}
where $\widehat{f}:V\oplus W\to L(V\oplus W,V\oplus W)$ is defined
by
\[
\widehat{f}(x,y)(\xi,\eta)=(\xi,f(y_{0}+y)\xi),\ \forall(x,y),(\xi,\eta)\in V\oplus W.
\]
\end{definition}

If the vector field $f\in C^{3}(W,L(V,W))$ satisfies the linear growth
and Lipschitz continuous conditions, then the existence and uniqueness
of a solution are ensured (see Theorem 6.2.1 and Corollary 6.2.2 in \cite{LQ02}). We will use $\Phi(y_{0},\bfX)$ to denote the unique solution
$\bfY$, call the map $\bfX\to\Phi(y_{0},\bfX)$ It\^o map defined by differential
equation (\ref{RDEXY}), whose Lipschitz continuity in $p$-variation
topology can be proved under the same conditions above (see Theorem 6.2.2 in \cite{LQ02}). This is an important result of It\^o maps in the framework
of differential equations driven by rough paths.

\subsubsection{Relation between differential equations driven by different rough
paths}

In this subsection, our main goal is to show the relationship of differential
equations driven by Stratonovich fractional Brownian rough path and
It\^o fractional Brownian rough path, respectively. Define $\Phi(x,\bfB^{\text{Str}})$
as the It\^o map of the differential equation
\begin{equation}
d\bfX=f(X)d\bfB^{\text{Str}},\ X_{0}=x,
\end{equation}
where $\bfB^{\text{Str}}$ is the Stratonovich fractional Brownian rough path,
and we use $d\bfB^{\text{Str}}$ to denote the equation driven by Stratonovich
rough path, which sometime we use $\circ d\bfB^{\text{Str}}$ instead. Respectively,
let $I(x,\bfB^{\text{It\^o}})$ denote the It\^o map of the differential equations driven
by the It\^o fractional rough path
\begin{equation}
d\bfX=f(X)d\bfB^{\text{It\^o}},\ X_{0}=x.\label{Ito-RDE}
\end{equation}
Now we want to ask what is the relationship between $I(x,\bfB^{\text{It\^o}})$ and
$\Phi(x,\bfB^{\text{Str}})$ or if the It\^o differential equation has a representation
in terms of a Stratonovich differential equation. First, we introduce
a geometric rough path $\bfB^{\text{Str},\varphi}$ defined by
\[
B^{\text{Str},\varphi}_{s,t}:=(B_{s,t},t^{2H}-s^{2H}),
\]
and
\[
\begin{split}\bbB^{\text{Str},\varphi}_{s,t} & :=\left(\bbB^{\text{Str}}_{s,t},\int_{s}^{t}B_{s,u}du^{2H}, \int_{s}^{t}(u^{2H}-s^{2H})dB_{u},\frac{1}{2}(t^{2H}-s^{2H})^{2}\right),
\end{split}
\]
where the cross integrals are Young integrals. Let $\Phi_{\widetilde{f}}(x,\bfB^{\text{Str},\varphi})$
be the It\^o map defined by the differential equation
\begin{equation}
d\bfX=\widetilde{f}(X)d\bfB^{\text{Str},\varphi},\ X_{0}=x,\label{Strat-RDE2}
\end{equation}
where $\widetilde{f}:\mathbb{R}^{e}\to L(\mathbb{R}^{d}\oplus\mathbb{R},\mathbb{R}^{e})$,
\[
\widetilde{f}(x)(\xi,\eta):=f(x)\xi-\frac{1}{2}\eta Df(x)f(x),
\]
for all $x\in\mathbb{R}^{e}$, $(\xi,\eta)\in\mathbb{R}^{d}\oplus\mathbb{R}$.
Namely, $\Phi_{\widetilde{f}}(\cdot,\bfB^{\text{Str},\varphi})$ is the It\^o map
of the rough differential equation
\begin{equation}
d\bfX=f(X)d\bfB^{\text{Str}}-\frac{1}{2}Df(X_{t})f(X_{t})dt^{2H}.\label{Strat-RDE22}
\end{equation}

\begin{theorem}\label{SI RDE} Let $f:\mathbb{R}^{e}\to L(\mathbb{R}^{d},\mathbb{R}^{e})$
be a $C^{4}$ vector field, and $I(x,\bfB^{\text{It\^o}})$, $\Phi_{\widetilde{f}}(x,\bfB^{\text{Str},\varphi})$
be the It\^o map defined by (\ref{Ito-RDE}), (\ref{Strat-RDE2}), respectively.
Then
\begin{equation}
I(x,\bfB^{\text{It\^o}})_{s,t}^{1}=\Phi_{\widetilde{f}}(x,\bfB^{\text{Str},\varphi})_{s,t}^{1}\label{RDE-trlt-1}
\end{equation}
and
\begin{equation}
I(x,\bfB^{\text{It\^o}})_{s,t}^{2}=\Phi_{\widetilde{f}}(x,\bfB^{\text{Str},\varphi})_{s,t}^{2}-\frac{1}{2}\int_{s}^{t}f(X_{u})\otimes f(X_{u})du^{2H},\label{RDE-trlt-2}
\end{equation}
where $X_{t}:=\Phi_{\widetilde{f}}(x,\bfB^{\text{Str},\varphi})_{0,t}^{1}$. \end{theorem}

The proof can be found in Appendix~\ref{appendix_SI_RDE}. As a remark, we can see that the differential equations (\ref{Ito-RDE})
and (\ref{Strat-RDE2}) (or (\ref{Strat-RDE22})) are not equivalent.
As far as the first level concerned, they define the same solution.
This agrees with classical stochastic differential equations (SDE)
driven by standard Brownian motion, i.e. we translate SDE (\ref{Ito-RDE})
into SDE (\ref{Strat-RDE22}) when $H=1/2$. However, in terms of
the second level, the two differential equations are different. The reason
is that $I(x,\bfB^{\text{It\^o}})^{2}$ can be viewed as the iterated integral of the
first level in It\^o's sense, while $\Phi_{\widetilde{f}}(x,\bfB^{\text{Str},\varphi})^{2}$
is the iterated integral of the first level in Stratonovich sense.

Besides, the relationship between $\Phi(x,\bfB^{\text{Str}})$ and $\Phi_{\widetilde{f}}(x,\bfB^{\text{Str},\varphi})$
is obvious. They are all understood in the Stratonovich sense, the
difference is their drift terms. By the continuity of It\^o's maps, they
can all be approximated in variation topology by the solution of differential
equations driven by piece-wise linear approximations of fBM and its
iterated path integrals. In summary, we can establish the relationship
between $I(x,\bfB^{\text{It\^o}})$, $\Phi(x,\bfB^{\text{Str}})$ and $\Phi_{\widetilde{f}}(x,\bfB^{\text{Str},\varphi})$.

\subsection{Examples and applications}

\label{example} In this subsection, we will give some interesting
examples as applications of fractional Brownian rough paths.

\noindent \textbf{Example 2.} Consider the differential equation in
dimension $d=1$,
\[
d\bfX=\sigma X\circ d\bfB^{\text{Str}}.
\]
Then $\bfX=(1,X,\bbX)$ is the solution of this differential equation,
where $X_{t}=\exp(\sigma B_{t})$,
\[
X_{s,t}=X_{t}-X_{s}=\exp(\sigma B_{t})-\exp(\sigma B_{s}),\ \ \
\]
\[
\bbX_{s,t}=\frac{1}{2}(X_{s,t})^{2}=\frac{1}{2}(\exp(\sigma B_{t})-\exp(\sigma B_{s}))^{2}.
\]

\noindent \textbf{Example 3.} Now we consider the It\^o rough differential
equation
\[
d\bfX=\sigma Xd\bfB^{\text{It\^o}}.
\]
(i) Set $X_{t}=\exp(\sigma B_{t}-\frac{1}{2}\sigma^{2}t^{2H})=:F(B_{t},t)$
and $X_{s,t}:=X_{t}-X_{s}$. Then
\[
X_{s,t}=\int_{s}^{t}\sigma X\rd\bfB^{\text{It\^o}},
\]
and
\[
\begin{split}\bbX_{s,t} & :=\int_{s}^{t}\sigma F(B,u)\circ \rrd\bfB^{\text{Str}}-\frac{\sigma^{2}}{2}\int_{s}^{t}(F(B_{u},u))^{2}du^{2H}\\
 & \ \ \ -\frac{\sigma^{2}}{2}\int_{s}^{t}\left(\int_{s}^{u}F(B_{r},r)dr^{2H}\right)dY_{u}-\frac{\sigma^{2}}{2}\int_{s}^{t}Y_{s,u}F(B_{u},u)du^{2H}\\
 & \ \ \ +\frac{\sigma^{4}}{4}\int_{s}^{t}\left(\int_{s}^{u}F(B_{r},r)dr^{2H}\right)F(B_{u},u)du^{2H},
\end{split}
\]
where $Y_{s,t}=\int_{s}^{t}\sigma F(B,u)\circ \rd\bfB^{\text{Str}}$, 
and the last four integrals are Young integrals.

\noindent \textbf{Example 4.} \textit{Geometric fBM} (or \textit{fractional
Black-Scholes model}). First consider the fractional Black-Scholes
model in Stratonovich sense
\begin{equation}
d\bfX=\mu X_{t}dt+\sigma X\circ d\bfB^{\text{Str}},
\end{equation}
where $\mu,\sigma$ are constants. The solution $\bfX$ can be constructed
as above. Set
\begin{equation}
X_{t}=\exp\left(\sigma B_{t}+\mu t\right)=:F(B_{t},t),
\end{equation}
and
\[
X_{s,t}=X_{t}-X_{s},
\]
\[
\bbX_{s,t}=\lim_{n\to\infty}\int_{s}^{t}F(B_{u}^{n},u)dF(B_{u}^{n},u),
\]

\noindent where $B_{u}^{n}$ is the linear approximation of $B$,
i.e.
\[
B_{u}^{n}=B_{t_{\ell-1}}+\frac{B_{t_{\ell}}-B_{t_{\ell-1}}}{t_{\ell}-t_{\ell-1}}(u-t_{\ell-1})
\]
on interval $[t_{\ell-1},t_{\ell}]$ and $\{t_{\ell}^{n},\ell=0,1,\cdots,n\}$
is any partition of $[s,t]$. We can verify $\bfX=(1,X,\bbX)$ is
the solution.

In this situation, the corresponding fractional Black-Scholes market
has arbitrage. We change the Stratonovich integral into an It\^o integral,
i.e. we consider the fractional differential equation in It\^o's sense
\begin{equation}
d\bfX=\mu X_{t}dt+\sigma Xd\bfB^{\text{It\^o}}.\label{fBS}
\end{equation}

We demonstrate later that the corresponding It\^o fractional Black-Scholes
market is arbitrage free in a restricted sense.\\
 (i) Let
\begin{equation}
X_{t}=\exp\left(\sigma B_{t}+\mu t-\frac{1}{2}\sigma^{2}t^{2H}\right)=:F(B_{t},t).\label{g-fBm}
\end{equation}
Then $X_{s,t}:=X_{t}-X_{s}$ and
\[
X_{s,t}=\int_{s}^{t}\sigma X\rd\bfB^{\text{It\^o}}+\int_{s}^{t}\mu X_{u}du.
\]
By the relation between Stratonovich integrals and It\^o integrals in
time dependent case, we have
\[
\begin{split}RHS & =\int_{s}^{t}\sigma F(B,u)\rd\bfB^{\text{It\^o}}+\int_{s}^{t}\mu F(B_{u},u)du\\
 & =\int_{s}^{t}\sigma F(B,u)\circ \rd\bfB^{\text{Str}}-\frac{\sigma^{2}}{2}\int_{s}^{t}F(B_{u},u)du^{2H}+\int_{s}^{t}\mu F(B_{u},u)du\\
 & =\int_{s}^{t}D_{x}F(B,u)\circ \rd\bfB^{\text{Str}}+\int_{s}^{t}D_{u}F(B_{u},u)du\\
 & =\int_{s}^{t}DF(\widetilde{B})\circ \rd\widetilde{\bfB}^{\text{Str}}=F(\widetilde{B}_{t})-F(\widetilde{B}_{s})\\
 & =F(B_{t},t)-F(B_{s},s)=X_{t}-X_{s}=LHS.
\end{split}
\]
(ii) Now set
\[
\begin{split}\bbX_{s,t} & :=\bbZ_{s,t}+\int_{s}^{t}\mu Z_{s,u}F(B_{u},u)du+\int_{s}^{t}\left(\int_{s}^{u}\mu F(B_{r},r)dr\right)dZ_{u}\\
 & \ \ \ \ \ \ \ \ \ +\int_{s}^{t}\left(\int_{s}^{u}\mu F(B_{r},r)dr\right)\mu F(B_{u},u)du,
\end{split}
\]
where
\[
\begin{split}Z_{s,t} & =X_{s,t}-\int_{s}^{t}\mu F(B_{u},u)du,\end{split}
\]
\[
\begin{split}\bbZ_{s,t} & \ =\int_{s}^{t}\sigma F(B,u)\circ \rrd\bfB^{\text{Str}}-\frac{\sigma^{2}}{2}\int_{s}^{t}(F(B_{u},u))^{2}du^{2H}\\
 & \ \ \ -\frac{\sigma^{2}}{2}\int_{s}^{t}\left(\int_{s}^{u}F(B_{r},r)dr^{2H}\right)dY_{u}-\frac{\sigma^{2}}{2}\int_{s}^{t}Y_{s,u}F(B_{u},u)du^{2H}\\
 & \ \ \ +\frac{\sigma^{4}}{4}\int_{s}^{t}\left(\int_{s}^{u}F(B_{r},r)dr^{2H}\right)F(B_{u},u)du^{2H},
\end{split}
\]
and $Y_{s,t}=\int_{s}^{t}\sigma F(B,u)\circ \rd\bfB^{\text{Str}}$, 
and all the integrals except ones involving $\circ d\bfB^{\text{Str}}$ are Young
integrals, and the integral against $\circ d\bfB^{\text{Str}}$ is Stratonovich rough
integral which can be computed by linear approximations. Actually,
by the relation between Stratonovich rough integrals and It\^o rough
integrals, we have
\[
Z_{s,t}=\int_{s}^{t}\sigma F(B,u)\rd\bfB^{\text{It\^o}}\quad \textrm{ and }\quad \bbZ_{s,t}=\int_{s}^{t}\sigma F(B,u)\rrd\bfB^{\text{It\^o}}.
\]
If we define $f(x,t)(\xi,\tau):=\sigma F(x,t)\xi+\mu F(x,t)\tau$,
$\widetilde{\bfB}=(\bfB,t)$ the space-time rough path of $\bfB$, then $\bbX_{s,t}=\int_{s}^{t}f(\widetilde{B})\rrd\widetilde{\bfB}^{\text{It\^o}}.$
Combining (i) and (ii), we have verified that $\bfX_{s,t}=\int_{s}^{t}f(\widetilde{B})d\widetilde{\bfB}^{\text{It\^o}}.$
The right hand side indeed coincides with the right hand side of the
differential equation (\ref{fBS}). So we have constructed the solution
of the It\^o fractional Black-Scholes equation (\ref{fBS}).

\begin{remark}\label{rmk-eqvlt} As a remark, we have two ways to understand
the integrals on the right hand side of above differential equation
with drift (\ref{fBS}). On the one hand, we can define $f(x,t)(\xi,\tau):=\sigma F(x,t)\xi+\mu F(x,t)\tau$,
then $\bfX_{s,t}=\int_{s}^{t}f(\widetilde{B})d\widetilde{\bfB}^{\text{It\^o}}$ is well
defined. On the other hand, we can define $g(x,t)(\xi,\tau):=\sigma F(x,t)\xi$,
$h_{t}:=\int_{0}^{t}\mu F(B_{u},u)du$ and see $\int_{s}^{t}\sigma F(B,u)d\bfB^{\text{It\^o}}$
as $\int_{s}^{t}g(\widetilde{B})d\widetilde{\bfB}^{\text{It\^o}}$ (This integral is
well defined). Then view the right hand side of differential equation
(\ref{fBS}) as a perturbation of $\int g(\widetilde{B})d\widetilde{\bfB}^{\text{It\^o}}$
by $h$. We want to say that the two ways are consistent, they give
the same results, i.e.
\begin{equation}
\int_{s}^{t}f(\widetilde{B})\rd\widetilde{\bfB}^{\text{It\^o}}=\int_{s}^{t}g(\widetilde{B})\rd\widetilde{\bfB}^{\text{It\^o}}+\int_{s}^{t}\mu F(B_{u},u)du,\label{consistent-1}
\end{equation}
\begin{equation}
\begin{split}\int_{s}^{t}f(\widetilde{B})\rrd\widetilde{\bfB}^{\text{It\^o}} & =\int_{s}^{t}g(\widetilde{B})\rrd\widetilde{\bfB}^{\text{It\^o}}+\int_{s}^{t}\mu Z_{s,u}^{1}F(B_{u},u)du +\int_{s}^{t}\left(\int_{s}^{u}\mu F(B_{r},r)dr\right)dZ_{u}\\
 & \ \ \ +\int_{s}^{t}\left(\int_{s}^{u}\mu F(B_{r},r)dr\right)\mu F(B_{u},u)du,
\end{split}
\label{consistent-2}
\end{equation}
where $Z_{s,t}:=\int_{s}^{t}g(\widetilde{B})\rd\widetilde{\bfB}^{\text{It\^o}}$,
$Z_{t}:=Z_{0}+Z_{0,t}$, and the last three integrals are Young Integral.
A proof of Equation (\ref{consistent-1}) and (\ref{consistent-2}) can be found in Appendix~\ref{appendix_rmk_eqvlt}.
\end{remark}

\subsection{Chain rule of fBM}

\label{CR-section}

This subsection is the continuity of Example 4 which plays an important
role in the remainder of the paper. Let $\bfX$ be the solution of fractional
Black-Scholes equation (\ref{fBS}) and $G$ is a good function, then
the integral $\int G(X,t)d\bfX$ is well defined. Since $X_{t}$ has
the explicit representation (\ref{g-fBm}),
\[
\int\sigma G(X,t)Xd\bfB^{\text{It\^o}}+\int\mu G(X_{t},t)X_{t}dt
\]
is defined in terms of the rough path $\bfB^{\text{It\^o}}$, which can be defined as
$\int f(\widetilde{B})d\widetilde{\bfB}^{\text{It\^o}}$, where
\[
f(x,t)(\xi,\eta):=f^{1}(x,t)\xi+f^{2}(x,t)\eta,
\]
\[
f^{1}(x,t)=\sigma G(F(x,t),t)F(x,t),
\]
\[
f^{2}(x,t)=\mu G(F(x,t),t)F(x,t),
\]
and
\[
F(x,t)=\exp(\sigma x-\sigma^{2}t^{2H}/2+\mu t).
\]
Heuristically, the two integrals should equal. Our aim in this subsection
is to show that they are indeed the same in this case. This kind of
formula is usually called \textit{Chain Rule}. Namely, we want to
show the theorem below.

\begin{theorem}\label{thm_chain_rule}(Chain rule) Let $\bfB^{\text{It\^o}}$ be the It\^o fractional Brownian rough
path with $H\in(\frac{1}{3},\frac{1}{2}]$, and $\bfX$ be the geometric
fBM with parameters $\sigma$ and $\mu$, $G$ be a function smooth
enough. Then
\begin{equation}
\int G(X,t)d\bfX=\int\sigma G(X,t)Xd\bfB^{\text{It\^o}}+\int\mu G(X_{t},t)X_{t}dt,
\end{equation}
where the integrals are understood as above.
\end{theorem}

The proof can be found in Appendix \ref{appendix_chain_rule}. We mention that with the Stratonovich rough paths, the chain rule
still holds by the same argument as It\^o rough paths above.

\subsection{Zero mean property of It\^o Integrals}

Now, we must verify whether the first level of our Itô integrals, as defined above, exhibit zero mean. Specifically, we need to confirm that:
\[
\mathbb{E}\left[\int_{s}^{t}f(B) \, d\bfB^{\text{Itô}}\right] = 0.
\]

\subsubsection{One dimensional case}

We prove it in the case that dimension $d=1$. First, we suppose
that $f$ has first and second continuous derivatives, by It\^o  formula
proved above, we can show that
\begin{equation}
\mathbb{E}\left[\int_{s}^{t}f'(B)\circ \rd\bfB^{\text{Str}}\right]=\mathbb{E}\left[\frac{1}{2}\int_{s}^{t}f''(B_{r})dr^{2H}\right].\label{f_1}
\end{equation}
The computation is routine, so we omit the details.

For the general case, set $F(x)=\int_{-\infty}^{x}f(y)dy$, $F(-\infty)=0$,
so that $F'(x)=f(x)$. By (\ref{f_1}) we get
\begin{equation}
\mathbb{E}\left[\int_{s}^{t}F'(B)\rd\bfB^{\text{It\^o}}\right]=\mathbb{E}\left[\frac{1}{2}\int_{s}^{t}F''(B_{r})dr^{2H}\right].
\end{equation}
Thus the expectation of the first level of the It\^o integral
\begin{equation}
\mathbb{E}\left[\int_{s}^{t}f(B)\rd\bfB^{\text{It\^o}}\right]=0.\label{f_11}
\end{equation}

\subsubsection{High dimensional case}

Now we can prove that eqn (\ref{f_1}) still holds when dimension
$d\geq2$ and for any function $F:\mathbb{R}^{d}\to L(\mathbb{R}^{d},\mathbb{R}^{e})$.
Let $\bfX_{s,t}:=\int_{s}^{t}F(B)d\bfB^{\text{It\^o}}$. The $i$-th component of first
level of this It\^o integral is
\[
\bfX_{s,t}^{1;i}=\left[\sum_{j=1}^{d}\int_{s}^{t}F^{ij}(B)\rd\bfB^{\text{It\^o},(j)}\right].
\]
As a remark here, the integral $\int_{s}^{t}F^{ij}(B)d\bfB^{\textbf{\text{It\^o}},(j)}$ is
well-defined, which can be understood as $\int_{s}^{t}\widetilde{F}^{ij}(B)d\bfB^{\text{It\^o}}$,
where $\widetilde{F}^{ij}(x_{1},\cdots,x_{d})(\xi_{1},\cdots,\xi_{d})=F^{ij}(x_{1},\cdots,x_{d})\xi_{j}$.
Therefore, we have the equation $\mathbb{E}\left[\int_{s}^{t}F^{ij}(B)\rd\bfB^{\text{It\^o},(j)}\right]=0$,
which yields that
\begin{equation}
\mathbb{E}\left[\int_{s}^{t}F(B)\rd\bfB^{\text{It\^o}}\right]=0.\label{f_d}
\end{equation}

\subsubsection{Zero mean property of time-dependent functions}

In this subsection, we will show that for the time-dependent function
$F(x,t)$ we can still have the mean zero property, i.e.
\begin{equation}
\mathbb{E}\left[\int_{s}^{t}F(B,u)\rd\bfB^{\text{It\^o}}\right]=0.\label{zeromean-1}
\end{equation}

As the time independent case, we first show the one dimensional case,
then by conditional expectation technique we conclude the high dimensional
cases. By the It\^o  formula (\ref{time-dept-Ito-1}) and Remark \ref{rmk-eqvlt},
we know that it is equivalent to
\begin{equation}
\begin{split}F(B_{t},t)-F(B_{s},s) & =\int_{s}^{t}D_{x}F(B,u)\rd\bfB^{\text{It\^o}}+\int_{s}^{t}D_{u}F(B_{u},u)du +\frac{1}{2}\int_{s}^{t}D_{x}^{2}F(B_{u},u)du^{2H}.
\end{split}
\end{equation}
Then in order to prove the zero mean property, we should verify that
\[
\mathbb{E}\left(F(B_{t},t)-F(B_{s},s)\right)=\int_{s}^{t}\mathbb{E}[D_{u}F(B_{u},u)]du+\frac{1}{2}\int_{s}^{t}\mathbb{E}[D_{x}^{2}F(B_{u},u)]du^{2H}.
\]
For the one dimension case, the left-hand side above can be computed
as the following
\begin{equation}
\mathbb{E}\left(F(B_{t},t)-F(B_{s},s)\right)=\int_{\mathbb{R}}(F(t^{H}x,t)-F(s^{H}x,s))\varphi(x)dx,\label{zeromean-111}
\end{equation}
where $\varphi$ is the standard normal probability density function.\\
On the other hand,
\begin{equation}
\begin{split} & \ \ \ \ \ \int_{s}^{t}\mathbb{E}[D_{u}F(B_{u},u)]du =\int_{s}^{t}\left[\int_{\mathbb{R}}\partial_{2}F(u^{H}x,u)\varphi(x)dx\right]du\\
 & =\int_{\mathbb{R}}\left[(F(t^{H}x,t)-F(s^{H}x,s))-\int_{s}^{t}x\partial_{1}F(u^{H}x,u)du^{H}\right]\varphi(x)dx.
\end{split}
\label{zeromean-222}
\end{equation}
and
\begin{equation}
\begin{split} & \ \ \ \ \ \frac{1}{2}\int_{s}^{t}\mathbb{E}[D_{x}^{2}F(B_{u},u)]du^{2H} =\frac{1}{2}\int_{s}^{t}\left[\int_{\mathbb{R}}\partial_{1}^{2}F(u^{H}x,u)\varphi(x)dx\right]du^{2H}\\
 & =\int_{s}^{t}\left[\int_{\mathbb{R}}\partial_{1}F(u^{H}x,u)x\varphi(x)dx\right]du^{H} =\int_{\mathbb{R}}\left[\int_{s}^{t}\partial_{1}F(u^{H}x,u)du^{H}\right]x\varphi(x)dx.
\end{split}
\label{zeromean-333}
\end{equation}
Combining (\ref{zeromean-111}),(\ref{zeromean-222}) and (\ref{zeromean-333}),
we thus obtain that
\begin{equation}
\mathbb{E}\left[\int_{s}^{t}D_{x}F(B,u)\rd\bfB^{\text{It\^o}}\right]=0.
\end{equation}
By using $\widetilde{F}(x,t)=\int_{-\infty}^{x}F(y,t)dy$, we finally
get (\ref{zeromean-1}) as for one dimensional case of one form above.

Now we turn to the high dimensional case. Let
\[
\bfX_{s,t}^{1;i}:=\left[\sum_{j=1}^{d}\int_{s}^{t}F^{ij}(B,u)\rd\bfB^{\text{It\^o},(j)}\right],\ i=1,\cdots,d.
\]
We need to prove that $\mathbb{E}\left[\int_{s}^{t}F^{ij}(B,u)\rd\bfB^{\text{It\^o},(j)}\right]=0$.
Let $F^{ij}=:f$ for simplicity. Then
\[
\begin{split} & \ \ \ \ \mathbb{E}\left[\int_{s}^{t}f(B,u)\rd\bfB^{\text{It\^o},(j)}\right]\\
 & =\mathbb{E}\left[\mathbb{E}\left[\left(\int_{s}^{t}f(x_{1},\cdots,B^{(j)},\cdots,x_{d},u)\rd\bfB^{\text{It\^o},(j)}\right)\right]_{x_{i}=B^{(i)},i\neq j}\right]\\
 & =\mathbb{E}\left[\mathbb{E}\left[\left(\int_{s}^{t}\widetilde{f}(B^{(j)},u)\rd\bfB^{\text{It\^o},(j)}\right)\right]_{x_{i}=B^{(i)},i\neq j}\right]\\
 & =0.
\end{split}
\]
Hence $\mathbb{E}\left[\bfX_{s,t}^{1;i}\right]=0$ for $i=1,\cdots,d$.
Thus we have proved that the expectation of time-dependent function
is also zero, i.e. eqn (\ref{zeromean-1}) holds.

\section{General It\^o integration for fBM}\label{sec-general-integration}

In the previous section, we introduced the It\^o integral for one-forms of fBM. For integrands that are not one-forms, we can similarly extend the concept of It\^o integration to accommodate them by adopting a renormalization strategy that transitions from the Stratonovich integral to the It\^o interpretation.

\subsection{It\^o integration for fOU}\label{sec-fOU-integration}

In the work by Qian and Xu (2024), \cite{QX2024}, the authors introduced the It\^o integral for fractional Ornstein-Uhlenbeck (fOU) processes and applied this integral to the problem of parameter estimation for the fOU process.

The multi-dimensional Ornstein-Uhlenbeck (OU) processes, driven by fractional Brownian motions (fBM), are commonly known as fractional Ornstein-Uhlenbeck (fOU) processes. These processes are defined by the solution to the stochastic differential equation (SDE):
\[
dX_t = -\Gamma X_t dt + \Sigma dB_t^H,\ X_0 = x_0.
\]
In this equation, \( B^H \) represents a \( d \)-dimensional fBM with a Hurst parameter \( H \) within the interval (0,1). The matrix \( \Gamma \) in \( \mathbb{R}^{d \times d} \) is the symmetric and positive-definite drift matrix, while \( \Sigma \) in \( \mathbb{R}^{d \times d} \) is the non-degenerate volatility matrix. This SDE is interpreted as a stochastic integral equation:
\[
X_t = x_0 - \int_0^t \Gamma X_s \, ds + \Sigma B_t^H,
\]
which admits a unique solution given by:
\[
X_t = e^{-\Gamma t} x_0 + \int_0^t e^{-\Gamma(t-s)} \Sigma \, dB_s^H.
\]
The integral on the right-hand side is interpreted using Young's integral approach. As a result, like the standard OU processes, the \( (X_t) \) series forms a Gaussian process.

To establish a non-geometric It\^o rough path enhancement compatible with fBM for analyzing fOU processes where \( \frac{1}{3} < H \leq \frac{1}{2} \), we define \( \varphi(t) := \frac{1}{2} It^{2H} - U(t) \), where:
\[
U(t) = H\Gamma\int_{0}^{t}\int_{0}^{s} e^{-\Gamma(s-u)} (s^{2H-1} - (s-u)^{2H-1}) \, du \, ds.
\]
This function demonstrates finite \( q \)-variation, with \( q = \frac{1}{2H} \), facilitating the definition of a non-geometric It\^o-type fractional Brownian rough path lift for \( B^H \) as:
\[
\mathbf{B}_{s,t}^{H,\text{It\^o}} = (1, B_{s,t}^{H}, \mathbb{B}_{s,t}^{H,\text{It\^o}}) := (1, B_{s,t}^{H}, \mathbb{B}_{s,t}^{H,\text{Str}} - \varphi_{s,t}),
\]
where \( \varphi_{s,t} = \varphi(t) - \varphi(s) \).

Based on the theory of differential equations driven by rough paths and the previously defined integrals, and assuming the coefficient matrices \(\Gamma\) and \(\Sigma\) commute for simplicity, we establish the following relationship:
\[
\int_{0}^{t}X_{s}d_{\mathfrak{R}_{1}}\mathbf{B}^{H,\text{It\^o}} = \int_{0}^{t}X_{s} \circ d_{\mathfrak{R}_{1}}\mathbf{B}^{H,\mathrm{Str}} - \Sigma\varphi(t).
\]
Given that
\[
\mathbb{E}\left[\int_{0}^{t}X_{s} \circ d_{\mathfrak{R}_{1}}\mathbf{B}^{H,\mathrm{Str}}\right] = \Sigma\varphi(t),
\]
we observe the zero expectation property for fOU processes, meaning:
\[
\mathbb{E}\left[\int_{0}^{t}X_{s}d_{\mathfrak{R}_{1}}\mathbf{B}^{H,\text{It\^o}}\right] = 0.
\]
This formulation reflects the expectation neutrality in the integrals of fOU processes under the specified rough path dynamics. For a more detailed exposition, we refer the reader to the work of Qian and Xu (2024), \cite{QX2024}.

\subsection{It\^o integration for RDE}\label{sec-RDE-integration}

In a more general context, consider a stochastic process \(\bfX\) that solves the rough differential equation (RDE) driven by a Stratonovich fractional Brownian rough path, formulated as:
\[
d\bfX = f(X) \circ d\mathbf{B}^{H,\text{Str}},\ X_{0} = x,
\]
where \(\mathbf{B}^{H,\text{Str}}\) is the Stratonovich fractional Brownian rough path with \( \frac{1}{3} < H \leq \frac{1}{2} \). To establish a corresponding It\^o integral for this process \(X\), we need to convert the Stratonovich integral into an It\^o integral, considering the nature of the driving fractional Brownian motion.

The transformation from a Stratonovich integral to an It\^o integral, in the context of fractional Brownian motion, is non-trivial. Similar to the above, we approach this by defining a non-geometric It\^o-type fractional Brownian rough path lift for \( B^H \) as follows:
\[
\mathbf{B}_{s,t}^{H,\text{It\^o}} = (1, B_{s,t}^{H}, \mathbb{B}_{s,t}^{H,\text{It\^o}}) := (1, B_{s,t}^{H}, \mathbb{B}_{s,t}^{H,\text{Str}} - \varphi_{s,t}),
\]
where \( \varphi_{s,t} = \varphi(t) - \varphi(s) \).

Consequently, the integral can be related to its Stratonovich counterpart:
\[
\int_{s}^{t}X_{u}d_{\mathfrak{R}_{1}}\mathbf{B}^{H,\text{It\^o}} = \int_{s}^{t}X_{u} \circ d_{\mathfrak{R}_{1}}\mathbf{B}^{H,\mathrm{Str}} - \int_{s}^{t} f(X_u)d \varphi(u).
\]
To ensure the expectation neutrality of the It\^o integral $\int_{s}^{t}X_{u}d_{\mathfrak{R}_{1}}\mathbf{B}^{H,\text{It\^o}}$, the function $\varphi(t)$ must satisfy the following equation: 
\[
\mathbb{E}\left[\int_{s}^{t}X_{u} \circ d_{\mathfrak{R}_{1}}\mathbf{B}^{H,\mathrm{Str}}\right] = \mathbb{E}\left[\int_{s}^{t} f(X_u)d \varphi(u)\right].
\]
Here, $\varphi(t)=\varphi(t;H,f)$ depends on the vector field $f$ and Hurst parameter $H$, $\varphi(0)=0$, and some regularity assumption needed. 

From a computational perspective, the function \(\varphi(t)\) can be approximated using a deep neural network \(\varphi_{\theta}(t;H,f)\) with parameters \(\theta\) that are trainable. To implement this, simulate a batch of \(N\) path samples \(\{X^i\}_{i=1}^N\) under the Stratonovich framework. The objective is to minimize the loss function defined as:
\[
L(\theta) =  \sum_{[s,t]\subset [0,T]}\left\| \frac{1}{N} \sum_{i=1}^N \left(\int_{s}^{t}X_{u}^i \circ d_{\mathfrak{R}_{1}}\mathbf{B}^{H,\text{Str}} - \int_{s}^{t} f(X_u^i) d \varphi_{\theta}(u) \right)\right\|^2,
\]
where both integrals are approximated by their respective discrete sums, i.e. 
\[
\int_{s}^{t}X_{u} \circ d_{\mathfrak{R}_{1}}\mathbf{B}^{H,\text{Str}}  \approx \sum_{\ell}\left(X_{t_{\ell-1}}B^H_{t_{\ell-1},t_{\ell}}+f(X_{t_{\ell-1}})\mathbb{B}^{H,Str}_{t_{\ell-1},t_{\ell}}\right),
\]
and 
\[
\int_{s}^{t} f(X_u) d \varphi_{\theta}(u) \approx \sum_{\ell} f(X_{t_{\ell-1}})(\varphi_{\theta}(t_{\ell})-\varphi_{\theta}(t_{\ell-1})).
\]
This leaves opportunities for further research in the future.

\section{Application in the fractional Black-Scholes model}

\label{Finance-sec} 

In this section, we extend the study of Example 4 presented in Section \ref{example} by considering the It\^o fractional Black-Scholes model, denoted by fBS, and its associated market. Our goal is to demonstrate that the market is arbitrage-free when the class of trading strategies is restricted. To achieve this, we first provide an arbitrage strategy under the Stratonovich fractional Black-Scholes market.

\subsection{Arbitrage strategy in Stratonovich fBS market}

Since the Stratonovich integral does not have zero mean property,
therefore fractional Black-Scholes market based on Stratonovich integral
suggests the possibility of existence of arbitrage. Shiryayev gave
an arbitrage trading strategy in \cite{Shir98} under Stratonovich
fBS market but driven by fBM with Hurst parameter $H>1/2$. We adapt
this strategy in our case when $H\in(\frac{1}{3},\frac{1}{2}]$.

The market has a stock (the risky asset) $\bfX$ whose price process
is $X_{t}:=\bfX_{0,t}^{1}$ at time $t$. We assume that $\bfX$ satisfies
the differential equation driven by Stratonovich fractional Brownian
rough path with $H\in(\frac{1}{3},\frac{1}{2}]$ as the first part
of example 4, i.e.
\begin{equation}
d\bfX=\mu X_{t}dt+\sigma X\circ d\bfB^{\text{Str}},\ X_{0}=x,\ t\in[0,T].\label{St-appl}
\end{equation}
The solution $\bfX$ of this equation has been constructed in example
4 and $X_{t}:=X_{0,t}^{1}=xe^{\sigma B_{t}+\mu t}$. It is assumed
that there is a money market (the risk-less asset) $M$, that is,
an asset whose price at time $t$ is not subject to uncertainty. Namely,
the price process $M_{t}$ satisfies the following equation
\begin{equation}
dM_{t}=rM_{t}dt,\ M_{0}=1,\ t\in[0,T],\label{money}
\end{equation}
where $r>0$ is a constant, i.e. $M_{t}=e^{rt}$.

A portfolio $(\gamma_{t},\zeta_{t})$ gives the number of units $\gamma_{t},\zeta_{t}$
held at time $t$ in the money market and stock market, respectively.
The value process $V_{t}\in\mathbb{R}$ of the portfolio is given
by
\begin{equation}
V_{t}=\gamma_{t}M_{t}+\zeta_{t}X_{t}.\label{V_t}
\end{equation}
The portfolio is called \textit{self-financing} if
\begin{equation}
V_{t}=V_{0}+\int_{0}^{t}\gamma_{s}dM_{s}+\int_{0}^{t}\zeta\circ \rd\bfX.
\end{equation}
Note that the second integral on the right hand side is the first
level of the Stratonovich integral against rough path $X$ defined
in (\ref{St-appl}).

Now consider the following portfolio
\begin{align}
\gamma_{t} & =1-e^{2\sigma B_{t}+2(\mu-r)t},\label{St-strategy-1}\\
\zeta_{t} & =2x^{-1}(e^{\sigma B_{t}+(\mu-r)t}-1),\label{St-strategy-2}
\end{align}
we will show that this trading strategy is an arbitrage one. First,
by (\ref{St-strategy-1}) and (\ref{St-strategy-2}), we get the value
process of the portfolio
\[
\begin{split}V_{t} & =(1-e^{2\sigma B_{t}+2(\mu-r)t})e^{rt}+2(e^{\sigma B_{t}+(\mu-r)t}-1)e^{\sigma B_{t}+\mu t}\\
 & =e^{rt}\left(e^{\sigma B_{t}+(\mu-r)t}-1\right)^{2}\geq0.
\end{split}
\]
By applying the basic principle/It\^o formula for Stratonovich integral
to $V_{t}=:f(B_{t},t)$, we have
\[
\begin{split}V_{t} & =V_{0}+\int_{0}^{t}re^{rs}\left(e^{\sigma B_{s}+(\mu-r)s}-1\right)^{2}ds\\
 & \ \ \ +\int_{0}^{t}2(\mu-r)e^{\sigma B_{s}+(\mu-r)s}e^{rs}\left(e^{\sigma B_{s}+(\mu-r)s}-1\right)ds\\
 & \ \ \ +\int_{0}^{t}2\sigma e^{rs}\left(e^{\sigma B_{s}+(\mu-r)s}-1\right)e^{\sigma B_{s}+(\mu-r)s}\circ \rd\bfB^{\text{Str}}_{s}\\
 & =\int_{0}^{t}r\gamma_{s}e^{rs}ds+\int_{0}^{t}\mu\zeta_{s}X_{s}ds+\int_{0}^{t}\sigma\zeta_{s}X_{s}\circ \rd\bfB^{\text{Str}}_{s}\\
 & =\int_{0}^{t}\gamma_{s}dM_{s}+\int_{0}^{t}\zeta_{s}\circ \rd\bfX_{s}.
\end{split}
\]
The last equality is by the chain rule of Stratonovich integral in
section \ref{CR-section}.

Hence, the portfolio (\ref{St-strategy-1}), (\ref{St-strategy-2})
is self-financing in this financial market. Note that the initial
payment at $t=0$ is $V_{0}=0$, but after that the value of this
portfolio is positive almost surely. This means one gets free lunch
with no risk.

\subsection{Arbitrage free under a class of trading strategies}

\label{arbitrage} Now we consider the It\^o fractional Black-Scholes
market. As for the Stratonovich fBS market, we suppose that the market
has a stock $\bfX$ (the risky asset) whose price process is $X_{t}:=\bfX_{0,t}^{1}$
but now it satisfies the differential equation driven by It\^o fractional
Brownian rough path.
\begin{equation}
d\bfX=\mu X_{t}dt+\sigma Xd\bfB^{\text{It\^o}},\ X_{0}=x,\ t\in[0,T].\label{Ito-fBS market}
\end{equation}
The solution is also constructed in example 4. The risk-less asset
money market $M$ satisfies the equation (\ref{money}), i.e. $M_{t}=e^{rt}$.

Suppose a portfolio $(\gamma_{t},\zeta_{t})$ gives the value process
$V_{t}\in\mathbb{R}$ by
\begin{equation}
V_{t}=\gamma_{t}M_{t}+\zeta_{t}X_{t}.\label{V_t}
\end{equation}
In this It\^o fBS market, we restrict the class of trading strategies.
We call a portfolio is \textit{admissible} if $\gamma_{t}=\gamma(X_{t},t)$,
and $\zeta_{t}=\zeta(X_{t},t)$. Besides, a portfolio is called \textit{self-financing}
if
\[
V_{t}=V_{0}+\int_{0}^{t}\gamma_{s}dM_{s}+\int_{0}^{t}\zeta \rd\bfX,
\]
where the second integral on the right hand side is the first level
of the It\^o integral against rough path $\bfX$ defined in (\ref{Ito-fBS market}).

Then by the chain rule of It\^o fractional Brownian rough path, we have
\[
\begin{split}V_{t} & =V_{0}+\int_{0}^{t}\gamma_{s}dM_{s}+\int_{0}^{t}\zeta \rd\bfX\\
 & =V_{0}+\int_{0}^{t}re^{rs}\gamma_{s}ds +\int_{0}^{t}\mu\zeta_{s}X_{s}ds+\int_{0}^{t}\sigma\zeta_{s}X_{s}\rd\bfB^{\text{It\^o}}
\end{split}
\]
By (\ref{V_t}) we also have
\begin{equation}
\gamma_{t}=e^{-rt}(V_{t}-\zeta_{t}X_{t}),
\end{equation}
plugging it into last equality, we get
\begin{equation}
V_{t}=V_{0}+\int_{0}^{t}rV_{s}ds+\int_{0}^{t}\sigma\zeta_{s}X_{s}\left(\frac{\mu-r}{\sigma}ds+\rd\bfB^{\text{It\^o}}\right).\label{V_t-2}
\end{equation}

In order to prove there is no arbitrage in this case, we first introduce
the Girsanov theorem for fBM with Hurst parameter $H\leq1/2.$

\subsection{Girsanov's theorem}

\label{Girsanov-sec}

The following version of Girsanov theorem for the fBM has been obtained
in (\cite{DU99}, Theorem 4.9), and we also suggest reader to see
\cite{BO03}, Theorem 4.1 and proof therein. In our case, we would
like to show that there is a new probability measure $\widehat{\mathbb{P}}$
such that
\begin{equation}
\widehat{B}_{t}=B_{t}+\frac{\mu-r}{\sigma}t,
\end{equation}
which is still an fBM under this measure $\widehat{\mathbb{P}}$.
This is what Girsanov theorem says in usual. Now let $K_{H}(t,s)$
be a square integrable kernel given by
\begin{equation}
K_{H}(t,s)=C_{H}\left[\frac{2}{2H-1}\left(\frac{t(t-s)}{s}\right)^{H-\frac{1}{2}}-\int_{s}^{t}\left(\frac{u(u-s)}{s}\right)^{H-\frac{1}{2}}\frac{du}{u}\right]1_{(0,t)}(s).
\end{equation}
Define the operator $K_{H}$ on $L^{2}([0,T])$ associated with the
kernel $K_{H}(t,s)$ as
\begin{equation}
(K_{H}f)(s)=\int_{0}^{T}f(t)K_{H}(t,s)dt.
\end{equation}

Given an adapted and integrable process $u=\{u_{t},\ t\in[0,T]\}$,
consider the transformation
\begin{equation}
\widehat{B}_{t}=B_{t}+\int_{0}^{t}u_{s}ds,\label{B-hat}
\end{equation}
since fBM $B$ can be represented by the integral along standard Brownian
motion $W$, we can write (\ref{B-hat}) into
\[
\widehat{B}_{t}=B_{t}+\int_{0}^{t}u_{s}ds=\int_{0}^{t}K_{H}(t,s)dW_{s}+\int_{0}^{t}u_{s}ds=\int_{0}^{t}K_{H}(t,s)d\widetilde{W}_{s},
\]
where $W_{t}$ is a standard Brownian motion and
\begin{equation}
\widetilde{W}_{t}=W_{t}+\int_{0}^{t}K_{H}^{-1}\left(\int_{0}^{\cdot}u_{r}dr\right)(s)ds.\label{B-tilde}
\end{equation}
By the standard Girsanov theorem for Brownian motion applied to (\ref{B-tilde}),
as a consequence, we have the following version of the Girsanov theorem
for the fBM with Hurst parameter $H\leq\frac{1}{2}$, which has obtained
in \cite{DU99}, \cite{NO02} and \cite{BO03}.

\begin{theorem}{(Girsanov theorem for fBM with $H\leq\frac{1}{2}$)}{(\cite{DU99},
Theorem 4.9; \cite{NO02}, Theorem 2; \cite{BO03}, Theorem 4.1)}\label{Girsanov} Let $B$
be a fBM with Hurst parameter $H\in(0,\frac{1}{2}]$, and
\[
v(s):=K_{H}^{-1}\left(\int_{0}^{\cdot}u_{r}dr\right)(s).
\]
Consider the shifted process (\ref{B-hat}). Assume that

(i) $\int_{0}^{T}u_{t}^{2}dt<\infty$, almost surely.

(ii) $\mathbb{E}(Z_{T})=1$, where
\[
Z_{T}=\exp\left(-\int_{0}^{T}v(s)dW_{s}-\frac{1}{2}\int_{0}^{T}\left(v(s)\right)^{2}ds\right),
\]
Then the shifted process $\widehat{B}$ is an $\mathcal{F}_{t}^{B}$-fBM
with Hurst parameter $H$ under the new probability measure $\widehat{\mathbb{P}}$
defined by $\frac{d\widehat{\mathbb{P}}}{d\mathbb{P}}=Z_{T}.$
\end{theorem}

\begin{remark} Here when $u$ satisfies the condition (i) in Theorem
\ref{Girsanov} with $H\leq\frac{1}{2}$, then $v=K_{H}^{-1}\left(\int_{0}^{\cdot}u_{r}dr\right)$
is well-defined, and $K_{H}^{-1}\left(\int_{0}^{\cdot}u_{r}dr\right)\in L^{2}([0,T])$,
where $K_{H}^{-1}$ is the inverse of the operator $K_{H}$. \end{remark}

In our case,
\[
\widehat{B}_{t}=B_{t}+\frac{\mu-r}{\sigma}t,
\]
so we can apply the Girsanov Theorem \ref{Girsanov} to it. Let $\mathbb{P}$
be the distribution of fBM $B$, and $\widehat{\mathbb{P}}$ be the
distribution constructed from $\mathbb{P}$ by Girsanov theorem. In
terms of $\widehat{B}_{t}$, we can write (\ref{V_t-2}), under $\widehat{\mathbb{P}}$,
as
\begin{equation}
V_{t}=V_{0}+\int_{0}^{t}rV_{s}ds+\int_{0}^{t}\sigma\zeta_{s}X\rd\widehat{\bfB}^{\text{It\^o}}.
\end{equation}
Transforming this equation, we have
\begin{equation}
e^{-rt}V_{t}=V_{0}+\sigma\int_{0}^{t}e^{-rs}\zeta_{s}X\rd\widehat{\bfB}^{\text{It\^o}},\ t\in[0,T].
\end{equation}
Taking expectation under the measure $\widehat{\mathbb{P}}$, we have
\begin{equation}
e^{-rT}\mathbb{E}_{\widehat{\mathbb{P}}}[V_{T}]=V_{0}+\mathbb{E}_{\widehat{\mathbb{P}}}\left[\sigma\int_{0}^{T}e^{-rs}\zeta_{s}X\rd\widehat{\bfB}^{\text{It\^o}}\right].
\end{equation}
By Girsanov's theorem and our zero mean property of integrals, we
get that
\begin{equation}
\mathbb{E}_{\widehat{\mathbb{P}}}\left[\int_{0}^{T}e^{-rs}\zeta_{s}X\rd\widehat{\bfB}^{\text{It\^o}}\right]=0.\label{expectation-0}
\end{equation}
Thus we may conclude that
\begin{equation}
e^{-rT}\mathbb{E}_{\widehat{\mathbb{P}}}\left[V_{T}\right]=V_{0}.
\end{equation}

Hence the probability measure $\widehat{\mathbb{P}}$ defined in Theorem
\ref{Girsanov} is a risk neutral measure. Then this It\^o fBS market
has no arbitrage under the class of admissible trading strategy. 

\subsection{Option pricing formula}

Furthermore, we aim to derive a pricing formula for the financial derivative $F$ at time $t=0$ under the risk-neutral measure $\widehat{\mathbb{P}}$. The arbitrage-free property of our It\^o fBS market under $\widehat{\mathbb{P}}$ is established in Theorem \ref{Girsanov}, which allows us to compute the price using the formula
\begin{equation}
V_{0}=e^{-rT}\mathbb{E}_{\widehat{\mathbb{P}}}\left[V_{T}\right].
\end{equation}
Once we have determined the price, we can state the following theorem.

\begin{theorem}{(Fractional Black-Scholes pricing formula)}\label{thm-bs-formula} The price
of claim $F(X_{T})$ under fractional Black-Scholes model and risk-neutral
measure $\widehat{\mathbb{P}}$ in Theorem \ref{Girsanov} is
\begin{equation}
V_{0}=e^{-rT}\int_{\mathbb{R}}F(X_{0}e^{\sigma T^{H}y+rT-\frac{1}{2}\sigma^{2}T^{2H}})\varphi(y)dy,
\end{equation}
where $\varphi(x)=\frac{1}{\sqrt{2\pi}}e^{-\frac{x^{2}}{2}}$ is the
standard normal density function.
\end{theorem}


The proof of Theorem~\ref{thm-bs-formula} can be found in Appendix~\ref{appendix_thm_bs_formula}. 
For the \textit{European call option}, $F(X_{T})=(X_{T}-K)^{+}$,
where $K>0$ is the striking price, we have the following conclusion.

\begin{corollary}{(European call)} The pricing formula of the European
call option under fractional Black-Scholes model and risk-neutral
measure $\widehat{\mathbb{P}}$ in Theorem \ref{Girsanov} is
\begin{equation}
V_{0}=X_{0}\left(1-\Phi(c_{-})\right)-Ke^{-rT}\left(1-\Phi(c_{+})\right),
\end{equation}
where
\[
c_{-}=\frac{1}{\sigma T^{H}}\log\left(\frac{K}{X_{0}}\right)-\frac{r}{\sigma}T^{1-H}-\frac{1}{2}\sigma T^{H},
\]
\[
c_{+}=\frac{1}{\sigma T^{H}}\log\left(\frac{K}{X_{0}}\right)-\frac{r}{\sigma}T^{1-H}+\frac{1}{2}\sigma T^{H},
\]
and $\Phi(x)=\int_{-\infty}^{x}\varphi(y)dy$ is standard normal distribution
function. \end{corollary}

\begin{figure}[ht]
\begin{center}
\includegraphics[width=12cm, height=8cm]{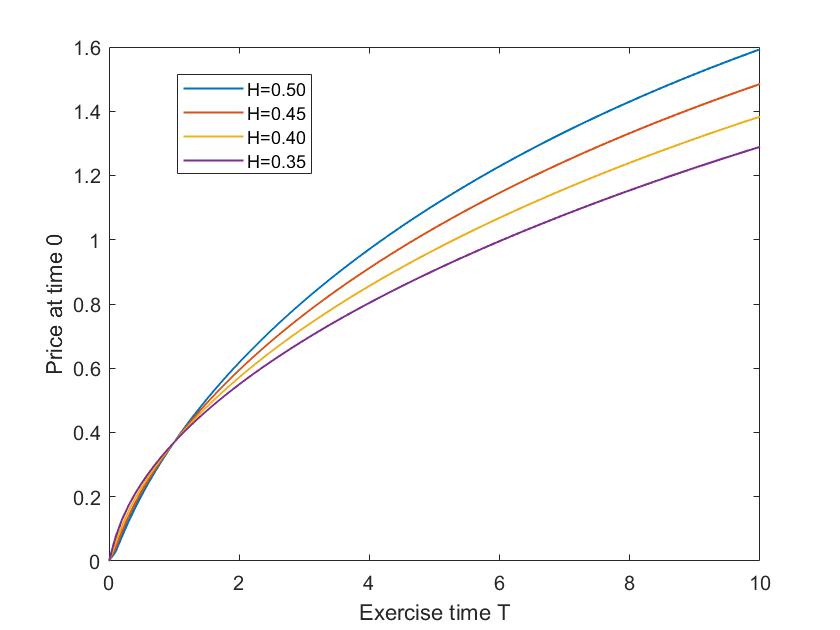}
\end{center}
\caption{Price of European call at time $0$ for different exercise time $T$ and varying Hurst parameter $H=0.50, 0.45, 0.40, 0.35$, where we take $\sigma=0.50$, $K=3$, $X_0=2.5$, $r=0.05$ as an example.}
\label{Price_fBS_fig}
\end{figure}


Let us comment on the pricing formula itself, which is similar
to that of the classical Black-Scholes model. The values are exactly
the same if the maturity time $T=1$, which seems strange at the first
glance. The Hurst parameter $H$ comes in to play a role only through
the exercise time $T$, and appears as a power in the exercise time
$T$ through $c_{\pm}$ and  to modify the volatility $\sigma^{2}T$
(for Black-Scholes' market) into $\sigma^{2}T^{2H}$ (for the fractional
BS market case), so that the intensity of the volatility is reduced
as $H<\frac{1}{2}$ due to the long time memory. In fact the numerics
$c_{\pm}$ can be rewritten as
\[
c_{\pm}=\frac{1}{\sigma T^{H}}\log\left(\frac{K}{X_{0}}\right)-\frac{rT}{\sigma T^{H}}\pm\frac{1}{2}\sigma T^{H}.
\]
Of course one has to understand the scale of the maturity $T$ in
time unit has no economic meaning, and its scale is in fact fixed
by the interest rate through $e^{rT}$, and therefore it looks natural
that $H$ should appear in the power of $T$ to have its effect on
the option pricing. See Figure \ref{Price_fBS_fig}, we show prices of European call option at time $0$ for different exercise time $T$ and Hurst parameter $H$.

\section{Application in the fractional Ornstein-Uhlenbeck process}\label{sec-fOU-estimate}

In Qian and Xu (2024), \cite{QX2024}, they introduced It\^o integration for fractional Ornstein-Uhlenbeck (fOU) processes, and applied the It\^o integration to the parameter estimation problem for the fOU process. In this section, we present the theoretical results regarding the parameter estimation problem for fOU processes. Furthermore, we showcase the efficacy of this rough path estimator via Monte Carlo simulations conducted on a 2-dimensional fOU process.

\subsection{Theoretical Results}
Let $X$ be a multi-dimensional fOU process, i.e. the solution to the following stochastic differential equation driven by multi-dimensional fractional Brownian motion
\begin{equation}
dX_{t}=-\Gamma X_{t}dt+\sigma dB_{t}^{H},\label{fOU-h-3}
\end{equation}
where $\Gamma$ is a parametric matrix, which is symmetric and positive-definite. We can construct an estimator based on either continuous or discrete observation via rough path theory. We suppose that the rough path enhancement $(X_{0,t}(\omega),\mathbb{X}_{0,t}(\omega))$ of the fOU process $X_{t}(\omega)$ can be observed continuously or discretely in It\^o sense. The main results are as follows. The proofs can be found in \cite{QX2024}.

\begin{theorem} 
(1) Suppose $\Gamma$ is a parametric
matrix and it is symmetric and positive-definite. Let $\widehat{\Gamma}_{t}$
be the rough path estimator of $\Gamma$: 
\begin{equation}\label{continuous-Gamma-estimator}
\widehat{\Gamma}_{t}^{T}\equiv -\left(\int_{0}^{t}X_{s}\otimes X_{s}ds\right)^{-1}\left(\int_{0}^{t}X_{s}\otimes d_{\mathfrak{R}_{1}}\mathbf{X}\right).
\end{equation}
Then
\begin{equation}
\widehat{\Gamma}_{t}\to\Gamma,\ a.s.,\ \text{as}\ t\to\infty.
\end{equation}
(2) Suppose we take the stochastic integral $\int_{0}^{t}X_{s}\otimes d_{\mathfrak{R}_{1}}\mathbf{X}$
in the rough path estimator $\widehat{\Gamma}_{t}$ 
as Stratonovich rough integral rather than It\^o rough integral as above,
we have
\begin{equation}
\widehat{\Gamma}_{t}\to0,\ a.s.,\ \text{as}\ t\to\infty.
\end{equation}
That is, we cannot use the Stratonovich rough integral to do this estimation problem. 
\end{theorem}

\begin{theorem}\label{thm-discrete-as} 
Suppose the fOU process $X$ with $H\in(\frac{1}{3},\frac{1}{2}]$ can be observed at discrete
time $\{t_{\ell}=\ell h,\ell=0,1,2,\cdots,n\}$ and as sample size
$n\to\infty$, $n$ and $h$ satisfy
\begin{equation}
nh\to\infty,\ h=h_{n}\to0,\ nh^{p}\to0,
\end{equation}
for some $p\in(1,\frac{1+H+\beta}{1+\beta})$, and $0<\beta<1$. Let
\begin{equation}\label{discrete-Gamma}
\widehat{\Gamma}_{n}^{T}\equiv-\left(\sum_{\ell=0}^{n}(X_{\ell h}\otimes X_{\ell h})h\right)^{-1}\left(\sum_{\ell=0}^{n-1}X_{\ell h}X_{\ell h,(\ell+1)h}+\mathbb{X}_{\ell h,(\ell+1)h}\right),
\end{equation}
where $\widehat{\Gamma}^{T}$ denotes transpose of matrix $\widehat{\Gamma}$.
Then
\begin{equation}
\widehat{\Gamma}_{n}\to\Gamma,\ a.s.
\end{equation}
as $n\to\infty$. 
\end{theorem}

As a remark, the explicit dependence of the L\'evy area of fOU processes on the drift parameter $\Gamma$ can complicate the application of the aforementioned rough path estimator in real-world observations. Nevertheless, as indicated by Equation (\ref{continuous-Gamma-estimator}) or (\ref{discrete-Gamma}), we can approach this issue by regarding it as an equation to be solved iteratively for $\Gamma$. This involves initializing with an approximate guess for the parameters, followed by a L\'evy area correction based on this preliminary estimate. Subsequently, we refine our parameter estimates in an iterative fashion.

\subsection{Monte Carlo Exercise}
In this subsection, we give a numerical example based on Monte Carlo simulation for a two dimensional
fOU process, i.e. the dynamic
\begin{equation}
dX_{t}=-\Gamma X_{t}dt+dB_{t}^{H},\ X_{0}=0,\ t\in[0,T],
\end{equation}
with parameter matrix
\[
\Gamma=\left(\begin{matrix}1 & 2\\
2 & 5
\end{matrix}\right).
\]
That is,
\begin{align*}
dX_{t}^{1} & =-(X_{t}^{1}+2X_{t}^{2})dt+dB_{t}^{H,1},\\
dX_{t}^{2} & =-(2X_{t}^{1}+5X_{t}^{2})dt+dB_{t}^{H,2}.
\end{align*}
We apply Euler scheme to draw $n$ equidistant samples $\{X_{h}(\omega),X_{2h}(\omega),\cdots,X_{nh}(\omega)\}$
on time horizon $T$ with frequency $h=\frac{T}{n}$.  In order to estimate the parameter matrix,
we should enhance sample paths to data in rough path sense. That is,
to get $\{(X_{\ell h,(\ell+1)h}(\omega),\mathbb{X}_{\ell h,(\ell+1)h}(\omega))_{\ell=0,1,\cdots,n-1}\}$
as our observation data for estimation.

\begin{table}[!ht]
\centering \caption{Mean and standard deviation of rough path estimators $\widehat{\Gamma}_{n}$ of the fOU process based on 1000 Monte Carlo simulation in dimension $d=2$, $n$ sample size, $T$ time horizon, $h$ sampling frequency, $H$ Hurst parameter, and true parameter matrix $\Gamma = \begin{pmatrix} 1 & 2 \\ 2 & 5 \end{pmatrix}$.}
\begin{tabular}{cccc}
\toprule
 & $n=2048$  & $n=4096$  & $n=8192$ \tabularnewline
\midrule
$H=0.50$  & $T=20,h=0.0098$  & $T=30,h=0.0073$  & $T=40,h=0.0049$ \tabularnewline
\midrule
\multirow{2}{*}{Mean}  & \multirow{2}{*}{$\left(\begin{matrix}1.0956 & 1.9720\\
1.9765 & 5.0381
\end{matrix}\right)$}  & \multirow{2}{*}{$\left(\begin{matrix}1.0548 & 1.9579\\
1.9918 & 5.0443
\end{matrix}\right)$}  & \multirow{2}{*}{$\left(\begin{matrix}1.0551 & 2.0091\\
1.9848 & 5.0146
\end{matrix}\right)$} \tabularnewline
 &  &  & \tabularnewline
\multirow{2}{*}{Std dev}  & \multirow{2}{*}{$\left(\begin{matrix}0.3350 & 0.6824\\
0.3374 & 0.7399
\end{matrix}\right)$}  & \multirow{2}{*}{$\left(\begin{matrix}0.2755 & 0.5955\\
0.2719 & 0.5590
\end{matrix}\right)$}  & \multirow{2}{*}{$\left(\begin{matrix}0.2429 & 0.5142\\
0.2313 & 0.4984
\end{matrix}\right)$} \tabularnewline
 &  &  & \tabularnewline
\midrule
$H=0.45$  & $T=13,h=0.0063$  & $T=22,h=0.0054$  & $T=40,h=0.0049$ \tabularnewline
\midrule
\multirow{2}{*}{Mean}  & \multirow{2}{*}{$\left(\begin{matrix}1.1163 & 1.9832\\
1.9889 & 5.1160
\end{matrix}\right)$}  & \multirow{2}{*}{$\left(\begin{matrix}1.0731 & 1.9943\\
1.9665 & 5.0106
\end{matrix}\right)$}  & \multirow{2}{*}{$\left(\begin{matrix}1.0449 & 2.0105\\
1.9744 & 4.9902
\end{matrix}\right)$} \tabularnewline
 &  &  & \tabularnewline
\multirow{2}{*}{Std dev}  & \multirow{2}{*}{$\left(\begin{matrix}0.4135 & 0.9030\\
0.4429 & 0.8453
\end{matrix}\right)$}  & \multirow{2}{*}{$\left(\begin{matrix}0.3172 & 0.6914\\
0.3034 & 0.6128
\end{matrix}\right)$}  & \multirow{2}{*}{$\left(\begin{matrix}0.2169 & 0.4743\\
0.2166 & 0.4517
\end{matrix}\right)$} \tabularnewline
 &  &  & \tabularnewline
\midrule
$H=0.40$  & $T=14,h=0.0068$  & $T=20,h=0.0049$  & $T=35,h=0.0043$ \tabularnewline
\midrule
\multirow{2}{*}{Mean}  & \multirow{2}{*}{$\left(\begin{matrix}1.1030 & 2.0077\\
1.9767 & 5.0416
\end{matrix}\right)$}  & \multirow{2}{*}{$\left(\begin{matrix}1.0585 & 1.9836\\
1.9895 & 5.0283
\end{matrix}\right)$}  & \multirow{2}{*}{$\left(\begin{matrix}1.0361 & 1.9946\\
1.9790 & 4.9894
\end{matrix}\right)$} \tabularnewline
 &  &  & \tabularnewline
\multirow{2}{*}{Std dev}  & \multirow{2}{*}{$\left(\begin{matrix}0.3745 & 0.8473\\
0.3878 & 0.7160
\end{matrix}\right)$}  & \multirow{2}{*}{$\left(\begin{matrix}0.3200 & 0.7466\\
0.3242 & 0.6152
\end{matrix}\right)$}  & \multirow{2}{*}{$\left(\begin{matrix}0.2275 & 0.5441\\
0.2367 & 0.4517
\end{matrix}\right)$} \tabularnewline
 &  &  & \tabularnewline
\midrule
$H=0.35$  & $T=14,h=0.0068$  & $T=20,h=0.0049$  & $T=30,h=0.0037$ \tabularnewline
\midrule
\multirow{2}{*}{Mean}  & \multirow{2}{*}{$\left(\begin{matrix}1.0724 & 1.9606\\
1.9588 & 4.9324
\end{matrix}\right)$}  & \multirow{2}{*}{$\left(\begin{matrix}1.0555 & 2.0065\\
1.9645 & 4.9839
\end{matrix}\right)$}  & \multirow{2}{*}{$\left(\begin{matrix}1.0298 & 1.9927\\
1.9796 & 4.9908
\end{matrix}\right)$} \tabularnewline
 &  &  & \tabularnewline
\multirow{2}{*}{Std dev}  & \multirow{2}{*}{$\left(\begin{matrix}0.3612 & 0.8816\\
0.3796 & 0.6448
\end{matrix}\right)$}  & \multirow{2}{*}{$\left(\begin{matrix}0.3074 & 0.7608\\
0.3234 & 0.5628
\end{matrix}\right)$}  & \multirow{2}{*}{$\left(\begin{matrix}0.2448 & 0.6067\\
0.2488 & 0.4453
\end{matrix}\right)$} \tabularnewline
 &  &  & \tabularnewline
\bottomrule
\end{tabular}\label{table_2d}
\end{table}

For dimension $d=2$, the continuous rough path estimator $\widehat{\Gamma}_{T}(\omega)=(\widehat{\Gamma}_{T}^{ij}(\omega))_{i,j=1,2}$
is given by
\begin{equation}
\begin{split}\widehat{\Gamma}_{T}^{ij}(\omega) & =-\frac{1}{V_{T}(X(\omega))}\left(\int_{0}^{T}(X_{s}^{3-j}(\omega))^{2}ds\int_{0}^{T}X_{s}^{j}(\omega)d_{\mathfrak{R}_{1}}\mathbf{X}^{i}(\omega)\right.\\
 & \ \ \ \ \ \ \ \ \ \ \ \ \ \ \ \left.-\int_{0}^{T}X_{s}^{i}(\omega)X_{s}^{3-i}(\omega)ds\int_{0}^{T}X_{s}^{3-j}(\omega)d_{\mathfrak{R}_{1}}\mathbf{X}^{i}(\omega)\right),
\end{split}
\end{equation}
where
\[
V_{T}(X(\omega))=\int_{0}^{T}(X_{s}^{1}(\omega))^{2}ds\int_{0}^{T}(X_{s}^{2}(\omega))^{2}ds-\left(\int_{0}^{T}X_{s}^{1}(\omega)X_{s}^{2}(\omega)ds\right)^{2},
\]
and all the rough integrals above are defined in our It\^o sense. Discretizing
every integral above, we obtain our high frequency rough path estimator.
The attention we need pay to is the cross term of rough integral,
i.e. L\'evy area $\int X^{i}d_{\mathfrak{R}_{1}}\mathbf{X}^{j}$ or
$\mathbb{X}^{ij}$. Since
\[
\mathbb{X}_{s,t}^{ij}=\mathbb{X}_{s,t}^{\text{Str},ij}-\varphi_{s,t},
\]
where $\varphi_{s,t}$ is defined in section \ref{sec-fOU-integration} and $\mathbb{X}_{s,t}^{\text{Str},ij}$
denotes the second level/L\'evy area of fOU rough path enhancement in
the Stratonovich sense. Now we can discretize the Stratonovich's L\'evy
area $\mathbb{X}^{\text{Str},ij}$ by trapezoidal scheme.
By this, we get $\{(X_{\ell h,(\ell+1)h}(\omega),\mathbb{X}_{\ell h,(\ell+1)h}(\omega))_{\ell=0,1,\cdots,n-1}\}$
as our discrete observation data for estimation.

We illustrate our two dimensional simulation results in Table \ref{table_2d}
below. In this case, we estimate the parameter matrix $\Gamma$ using
the simulated data $\{(X_{\ell h,(\ell+1)h}(\omega),\mathbb{X}_{\ell h,(\ell+1)h}(\omega))_{\ell=0,1,\cdots,n-1}\}$.
We draw 1000 sample paths by Monte Carlo simulation.

In Table \ref{table_2d}, every component of ``Mean" denotes average
of the value of respective estimator based on 1000 Monte Carlo simulation.
And the component of ``Standard deviation (Std dev)" represents the
fluctuation of estimation of parameter with corresponding index. One
could see that, under proper time horizon $T$, sample size $n$ and
frequency $h$, the rough path estimator performs very well and the
results are quite stable.


As a remark, in dimension $d=2$, the performance of discrete estimator becomes a little sensitive to sampling mode. One should adhere to the conditions about time horizon $T$, sample size $n$ and frequency $h$ in Theorem \ref{thm-discrete-as} in order to obtain better estimated values. In Table \ref{table_2d}, we set frequency $h$ becomes smaller as sample size $n$ becomes larger and Hurst parameter $H$ smaller.

\newpage 

\appendix

\setcounter{equation}{0}
\numberwithin{equation}{section}
\setcounter{figure}{0}
\numberwithin{figure}{section}

\renewcommand{\thesection}{Appendix A.}
\section{Proofs}
\renewcommand{\thesection}{A}

In this appendix, we provide the proofs for the results stated in the main text. For clarity, unless specified otherwise, we may use the notation $\bfX=(1, X, \bbX)$ interchangeably with $\bfX=(1, \bfX^{1}, \bfX^{2})$ for rough path $\bfX$, as long as this substitution does not lead to any confusion.

\subsection{Proof of Theorem~\ref{SI-ingtral-thm-notime} and \ref{SI-ingtral-thm}}\label{appendix_SI_correction}

\begin{proof} 
(i) By definition of our integral, we have
\[
\begin{split} & \ \ \ \ \int_{s}^{t}F(B,u)\rd\bfB^{\text{It\^o}}=\int_{s}^{t}f(\widetilde{B})\rd\widetilde{\bfB}^{\text{It\^o}}\\
 & =\lim_{|\mathcal{P}|\to0}\sum_{[u,v]\in\mathcal{P}}f(\widetilde{B}_{u})\widetilde{B}_{u,v}+Df(\widetilde{B}_{u})\widetilde{\bbB}_{u,v}^{\text{It\^o}}\\
 & =\lim_{|\mathcal{P}|\to0}\sum_{[u,v]\in\mathcal{P}}F(B_{u},u)B_{u,v}+D_{x}F(B_{u},u)\bbB_{u,v}^{\text{It\^o}} +D_{u}F(B_{u},u)\int_{u}^{v}B_{u,r}dr.
\end{split}
\]
Since
\[
\left|\int_{u}^{v}B_{r}dr-B_{u}(v-u)\right|=o(|v-u|)=o(|\mathcal{P}|),
\]
so that
\[
\int_{s}^{t}F(B,u)\rd\bfB^{\text{It\^o}}=\lim_{|\mathcal{P}|\to0}\sum_{[u,v]\in\mathcal{P}}F(B_{u},u)B_{u,v}+D_{x}F(B_{u},u)\bbB_{u,v}^{\text{It\^o}}.
\]
Therefore we may conclude that
\[
\begin{split} & \ \ \ \ \int_{s}^{t}F(B,u)\rd\bfB^{\text{It\^o}}\\
 & =\lim_{|\mathcal{P}|\to0}\sum_{[u,v]\in\mathcal{P}}F(B_{u},u)B_{u,v}+D_{x}F(B_{u},u)(\bbB^{\text{Str}}_{u,v}-\frac{1}{2}I(v^{2H}-u^{2H}))\\
 & =\lim_{|\mathcal{P}|\to0}\sum_{[u,v]\in\mathcal{P}}\left(F(B_{u},u)B_{u,v}+D_{x}F(B_{u},u)\bbB^{\text{Str}}_{u,v}\right)\\
 & \ \ \ \ \ \ \ \ \ -\frac{1}{2}\lim_{|\mathcal{P}|\to0}\sum_{[u,v]\in\mathcal{P}}D_{x}F(B_{u},u)(v^{2H}-u^{2H})\\
 & =\int_{s}^{t}F(B,u)\circ \rd\bfB^{\text{Str}}-\frac{1}{2}\int_{s}^{t}D_{x}F(B_{u},u)du^{2H}.
\end{split}
\]
(ii) Now we show the second relation (\ref{SI-22}). It follows in
the similar way as (i).
\[
\begin{split} & \ \ \ \ \int_{s}^{t}F(B,u)\rrd\bfB^{\text{It\^o}}=\int_{s}^{t}f(\widetilde{B})\rrd\widetilde{\bfB}^{\text{It\^o}}\\
 & =\lim_{|\mathcal{P}|\to0}\sum_{[u,v]\in\mathcal{P}}X^{\text{It\^o}}_{s,u}\otimes X^{\text{It\^o}}_{u,v}+f(\widetilde{B}_{u})\otimes f(\widetilde{B}_{u})\widetilde{\bbB}_{u,v}^{\text{It\^o}}.
\end{split}
\]
By eqn (\ref{SI-11}), and $f(\widetilde{B}_{u})\otimes f(\widetilde{B}_{u})\widetilde{\bbB}_{u,v}^{\text{It\^o}}=F(B_{u},u)\otimes F(B_{u},u)\bbB_{u,v}^{\text{It\^o}}$,
we obtain
\[
\begin{split} & \ \ \ \ \int_{s}^{t}F(B,u)\rrd\bfB^{\text{It\^o}}\\
 & =\lim_{|\mathcal{P}|\to0}\sum_{[u,v]\in\mathcal{P}}\left(X^{\text{Str}}_{s,u}-\frac{1}{2}\int_{s}^{u}D_{x}F(B_{r},r)dr^{2H}\right)\\
 & \ \ \ \ \ \ \ \otimes\left(X^{\text{Str}}_{u,v}-\frac{1}{2}\int_{u}^{v}D_{x}F(B_{r},r)dr^{2H}\right)\\
 & \ \ \ \ \ \ \ +F(B_{u},u)\otimes F(B_{u},u)\left(\bbB^{\text{Str}}_{u,v}-\frac{1}{2}(v^{2H}-u^{2H})\right)\\
 & =\lim_{|\mathcal{P}|\to0}\sum_{[u,v]\in\mathcal{P}}X^{\text{Str}}_{s,u}\otimes X^{\text{Str}}_{u,v}+F(B_{u},u)\otimes F(B_{u},u)\bbB^{\text{Str}}_{u,v}\\
 & \ \ \ \ \ \ \ -\frac{1}{2}F(B_{u},u)\otimes F(B_{u},u)\left(v^{2H}-u^{2H}\right)\\
 & \ \ \ \ \ \ \ -\frac{1}{2}\int_{s}^{u}D_{x}F(B_{r},r)dr^{2H}\otimes X^{\text{Str}}_{u,v}\\
 & \ \ \ \ \ \ \ -X^{\text{Str}}_{s,u}\otimes\frac{1}{2}\int_{u}^{v}D_{x}F(B_{r},r)dr^{2H}\\
 & \ \ \ \ \ \ \ +\frac{1}{4}\int_{s}^{u}D_{x}F(B_{r},r)dr^{2H}\int_{u}^{v}D_{x}F(B_{r},r)dr^{2H}.
\end{split}
\]
Since $B_{t}$ has finite $p$-variation with $p>1/H$, therefore
$\bbB$ is of finite $p/2$-variation, $F(B_{t},t)$, $D_{x}F(B_{t},t)$
have finite $p$-variations, $\int_{0}^{\cdot}D_{x}F(B_{r},r)dr^{2H}$
has finite $1/2H$-variation and $X^{\text{Str}}$ has finite $p$-variation.
Since $\frac{1}{H}<p<3$ and $\frac{1}{3}<H<\frac{1}{2}$, so that
$\frac{1}{p}+2H>1$, and the last four sums on the right hand side
converge to the Young integral. Eqn (\ref{SI-22}) therefore follows
immediately. 
\end{proof}

\subsection{Proof of Theorem~\ref{SI RDE}}\label{appendix_SI_RDE}

\begin{proof} 
Suppose $\{\overline{B}_{t},\ t\geq0\}$ is a piece-wise
linear/smooth approximation with finite variation of fractional Brownian
motion $\{B_{t},\ t\geq0\}$. Set $\overline{B}_{s,t}=\overline{B}_{t}-\overline{B}_{s}$
and let $\overline{\bbB}_{s,t}$ be the difference of the iterated
integral over $[s,t]$ of $\overline{B}$ and $\frac{1}{2}(t^{2H}-s^{2H})$.
Consider the rough differential equation
\begin{equation}
d\bfX=f(X)d\overline{\bfB},
\end{equation}
that is, the integral equation
\begin{equation}
\bfZ=\int\widehat{f}(Z)d\bfZ,\ \pi_{d}(\bfZ)=\overline{\bfB},
\end{equation}
where $\widehat{f}(x,y)(\xi,\eta):=(\xi,f(y)\xi)$, and $\pi_{d}$
is projection operator to $\mathbb{R}^{d}$, which is solved by the
Picard iteration, that is
\[
\bfZ(n+1)=\int\widehat{f}(Z(n))d\bfZ(n),\quad \bfZ(0)=(\overline{\bfB},0).
\]
More precisely, define almost rough paths
\begin{equation}
\widehat{\bfZ}(n+1)_{s,t}^{1}:=\widehat{f}(Z(n)_{s})\bfZ(n)_{s,t}^{1}+D\widehat{f}(Z(n)_{s})\bfZ(n)_{s,t}^{2},\label{arp-1}
\end{equation}
\begin{equation}
\widehat{\bfZ}(n+1)_{s,t}^{2}:=\widehat{f}(Z(n)_{s})\otimes\widehat{f}(Z(n)_{s})\bfZ(n)_{s,t}^{2},\label{arp-2}
\end{equation}
and define the corresponding rough paths
\begin{equation}
\bfZ(n+1)_{s,t}^{1}=\lim_{|\mathcal{P}|\to0}\sum_{[u,v]\in\mathcal{P}}\widehat{\bfZ}(n+1)_{u,v}^{1},\label{rp-Z1}
\end{equation}
\begin{equation}
\bfZ(n+1)_{s,t}^{2}=\lim_{|\mathcal{P}|\to0}\sum_{[u,v]\in\mathcal{P}}\bfZ(n)_{s,u}^{1}\otimes \bfZ(n)_{u,v}^{1}+\widehat{\bfZ}(n+1)_{u,v}^{2},\label{rp-Z2}
\end{equation}
where $\mathcal{P}$ is a partition of the interval $[s,t]$. Then
\begin{equation}
|\bfZ(n)_{s,t}^{i}-\widehat{\bfZ}(n)_{s,t}^{i}|\leq\omega(s,t)^{\theta},\ i=1,2,\ \theta>1,\label{control-ZZhat}
\end{equation}
for some control $\omega$. 

(i) Now we prove (\ref{RDE-trlt-1}). Set $\widehat{\bfZ}(n)_{s,t}^{1}=(\overline{\bfB}_{s,t}^{1},\bfX(n)_{s,t}^{1})$,
by definition of $\widehat{f}$ and (\ref{arp-1}), we have
\[
\begin{split}\widehat{\bfZ}(n+1)_{s,t}^{1} & =(\overline{\bfB}_{s,t}^{1},f(X(n)_{s})\overline{\bfB}_{s,t}^{1})+D\widehat{f}(Z(n)_{s})\bfZ(n)_{s,t}^{2}\\
 & \simeq(\overline{\bfB}_{s,t}^{1},f(X(n)_{s})\overline{\bfB}_{s,t}^{1})+D\widehat{f}(Z(n)_{s})\widehat{\bfZ}(n)_{s,t}^{2}\\
 & \simeq(\overline{\bfB}_{s,t}^{1},f(X(n)_{s})\overline{\bfB}_{s,t}^{1})+(0,Df(X(n)_{s})f(X(n)_{s})\overline{\bfB}_{s,t}^{2}\\
 & \simeq(\overline{\bfB}_{s,t}^{1},f(X(n)_{s})\overline{\bfB}_{s,t}^{1})-\frac{1}{2}(0,Df(X(n)_{s})f(X(n)_{s})(t^{2H}-s^{2H}),
\end{split}
\]
where $\simeq$ means the error can be controlled by $\omega(s,t)^{\theta}$
with $\theta>1$. Hence,
\[
\bfX(n+1)_{s,t}^{1}\simeq f(X(n)_{s})\overline{\bfB}_{s,t}^{1}-\frac{1}{2}Df(X(n)_{s})f(X(n)_{s})(t^{2H}-s^{2H}).
\]
Since $\overline{\bfB}^{1}$ has finite variation and $\sum_{[u,v]\in\mathcal{P}}\bfX(n)_{u,v}^{1}=\bfX(n)_{s,t}^{1}$,
the formula above implies that
\[
\bfX(n+1)_{s,t}^{1}=\lim_{|\mathcal{P}|\to0}\sum_{[u,v]\in\mathcal{P}}f(X(n)_{u})\overline{\bfB}_{u,v}^{1}-\frac{1}{2}Df(X(n)_{u})f(X(n)_{u})(v^{2H}-u^{2H}),
\]
as $n$ goes to infinity, the limit above is identified as
\[
X_{t}=X_{s}+\int_{s}^{t}f(X_{u})\rd\overline{\bfB}_{u}-\frac{1}{2}\int_{s}^{t}Df(X_{u})f(X_{u})du^{2H}.
\]
On the other hand,
\[
\lim_{n\to\infty}\bfX(n)_{s,t}^{1}=\Phi(\cdot,\overline{\bfB})_{s,t}^{1},
\]
so we get (\ref{RDE-trlt-1}) when the system is driven by $\overline{\bfB}$.
By the continuity of It\^o's maps, we conclude that (\ref{RDE-trlt-1})
holds in real fractional Brownian rough path case.

(ii) Similarly, by (\ref{rp-Z2}) and the continuity of It\^o's maps again, we conclude that (\ref{RDE-trlt-2}) holds. 
\end{proof}

\subsection{Proof of Theorem~\ref{thm_chain_rule}}\label{appendix_chain_rule}

\begin{proof} Let $\bfY=\int G(X,t)d\bfX$, $\mathbf{H}=\int f(\widetilde{B})d\widetilde{\bfB}$.
Then their associated almost rough path are $\widehat{\bfY}$, $\widehat{\mathbf H}$
respectively, where
\begin{align*}
\widehat{\bfY}_{s,t}^{1} & =G(X_{s},s)\bfX_{s,t}^{1}+D_{x}G(X_{s},s)\bfX_{s,t}^{2},\\
\widehat{\bfY}_{s,t}^{2} & =G(X_{s},s)\otimes G(X_{s},s)\bfX_{s,t}^{2},
\end{align*}
and
\begin{align*}
\widehat{\mathbf H}_{s,t}^{1} & =f^{1}(B_{s},s)\bfB_{s,t}^{1}+f^{2}(B_{s},s)(t-s)+D_{x}f^{1}(B_{s},s)\bfB_{s,t}^{2},\\
\widehat{\mathbf H}_{s,t}^{2} & =f^{1}(B_{s},s)\otimes f^{1}(B_{s},s)\bfB_{s,t}^{2}.
\end{align*}
We should prove that
\begin{equation}
|\widehat{\bfY}_{s,t}^{i}-\widehat{\mathbf H}_{s,t}^{i}|\leq\omega(s,t)^{\theta},\ i=1,2,\ \forall\ s<t,\ \exists\ \theta>1,\label{chain-YH}
\end{equation}
for some control $\omega$. If the difference of two quantities is
controlled by $\omega$ like (\ref{chain-YH}), we use the symbol
$\simeq$ to represent it. So we should show $\widehat{\bfY}_{s,t}^{i}\simeq\widehat{\mathbf H}_{s,t}^{i}$.

Denote $\bfZ=\int\sigma F(B,t)d\bfB$, $h_{s,t}=\int_{s}^{t}\mu F(B_{t},t)dt$,
and $\widehat{\bfZ}$, $\widehat{h}$ denote their respective almost
rough path, that is,
\begin{align*}
\widehat{\bfZ}_{s,t}^{1} & =\sigma F(B_{s},s)\bfB_{s,t}^{1}+\sigma^{2}F(B_{s},s)\bfB_{s,t}^{2},\\
\widehat{\bfZ}_{s,t}^{2} & =\sigma^{2}F(B_{s},s)\otimes F(B_{s},s)\bfB_{s,t}^{2},\\
\widehat{h}_{s,t} & =\mu F(B_{s},s)(t-s).
\end{align*}
It is easy to verify that the almost rough path associated with $\bfX$
is given by
\begin{equation}
(1,\bfZ_{s,t}^{1}+h_{s,t},\bfZ_{s,t}^{2})\simeq(1,\bfZ_{s,t}^{1}+\widehat{h}_{s,t},\bfZ_{s,t}^{2}).\label{arp-X}
\end{equation}
Two quantities on both sides in (\ref{arp-X}) are all almost rough
paths. So we have the following relations:
\begin{align}
\bfX_{s,t}^{1} & =\bfZ_{s,t}^{1}+h_{s,t}\simeq\widehat{\bfZ}_{s,t}^{1}+\widehat{h}_{s,t},\\
\bfX_{s,t}^{2} & \simeq \bfZ_{s,t}^{2}\simeq\widehat{\bfZ}_{s,t}^{2}.
\end{align}

(i) For the first level,
\[
\begin{split}\widehat{\bfY}_{s,t}^{1} & =G(X_{s},s)\bfX_{s,t}^{1}+D_{x}G(X_{s},s)\bfX_{s,t}^{2}\\
 & =G(F(B_{s},s),s)\bfX_{s,t}^{1}+\partial_{1}G(F(B_{s},s),s)\bfX_{s,t}^{2}\\
 & \simeq G(F(B_{s},s),s)(\widehat{\bfZ}_{s,t}^{1}+\widehat{h}_{s,t})+\partial_{1}G(F(B_{s},s),s)\widehat{\bfZ}_{s,t}^{2}\\
 & =G(F(B_{s},s),s)(\sigma F(B_{s},s)\bfB_{s,t}^{1}\\
 & \ \ \ +\sigma^{2}F(B_{s},s)\bfB_{s,t}^{2}+\mu F(B_{s},s)(t-s))\\
 & \ \ \ +\partial_{1}G(F(B_{s},s),s)(\sigma^{2}F(B_{s},s)\otimes F(B_{s},s)\bfB_{s,t}^{2}).
\end{split}
\]
Since
\[
D_{x}f^{1}(x,s)=\sigma^{2}\partial_{1}G(F(x,s),s)F(x,s)\otimes F(x,s)+\sigma^{2}G(F(x,s),s)F(x,s),
\]
so that
\[
\widehat{\bfY}_{s,t}^{1}\simeq f^{1}(B_{s},s)\bfB_{s,t}^{1}+f^{2}(B_{s},s)(t-s)+D_{x}f^{1}(B_{s},s)\bfB_{s,t}^{2}=\widehat{\mathbf H}_{s,t}^{1}.
\]
Thus we have proved the first part of the claim.

(ii) For the second level paths, we have
\[
\begin{split}\widehat{\bfY}_{s,t}^{2} & =G(X_{s},s)\otimes G(X_{s},s)\bfX_{s,t}^{2}\\
 & =G(F(B_{s},s),s)\otimes G(F(B_{s},s),s)\bfX_{s,t}^{2}\\
 & \simeq G(F(B_{s},s),s)\otimes G(F(B_{s},s),s)\widehat{\bfZ}_{s,t}^{2}\\
 & =G(F(B_{s},s),s)\otimes G(F(B_{s},s),s)(\sigma^{2}F(B_{s},s)\otimes F(B_{s},s)\bfB_{s,t}^{2})\\
 & =\sigma^{2}G(F(B_{s},s),s)F(B_{s},s)\otimes G(F(B_{s},s),s)F(B_{s},s)\bfB_{s,t}^{2}\\
 & =f^{1}(B_{s},s)\otimes f^{1}(B_{s},s)\bfB_{s,t}^{2}=\widehat{\mathbf H}_{s,t}^{2},
\end{split}
\]
which thus completes our proof.
\end{proof}

\subsection{Proof of Remark~\ref{rmk-eqvlt}}\label{appendix_rmk_eqvlt}

\begin{proof} (i) For the first level, we have
\[
\begin{split}\int_{s}^{t}f(\widetilde{B})\rd\widetilde{\bfB}^{\text{It\^o}} & =\lim_{|\mathcal{P}|\to0}\sum_{[u,v]\in\mathcal{P}}f(\widetilde{B}_{u})\widetilde{B}_{u,v}+Df(\widetilde{B}_{u})\widetilde{\bbB}_{u,v}^{\text{It\^o}}\\
 & =\lim_{|\mathcal{P}|\to0}\sum_{[u,v]\in\mathcal{P}}\sigma F(B_{u},u)B_{u,v}+\mu F(B_{u},u)(v-u)\\
 & \ \ \ \ +\sigma D_{x}F(B_{u},u)\bbB_{u,v}^{\text{It\^o}}+\sigma D_{u}F(B_{u},u)\int_{u}^{v}(B_{r}-B_{u})dr\\
 & \ \ \ \ +\mu D_{x}F(B_{u},u)\int_{u}^{v}(r-u)dB_{r}+\frac{\mu}{2}D_{u}F(B_{u},u)(v-u)^{2},
\end{split}
\]
Since
\[
\left|\int_{u}^{v}B_{r}dr-B_{u}(v-u)\right|=o(|v-u|)=o(|\mathcal{P}|)
\]
and
\[
\left|\int_{u}^{v}rdB_{r}-u(B_{v}-B_{u})\right|=o(|v-u|)=o(|\mathcal{P}|),
\]
we therefore have
\[
\begin{split}\int_{s}^{t}f(\widetilde{B})\rd\widetilde{\bfB}^{\text{It\^o}} & =\lim_{|\mathcal{P}|\to0}\sum_{[u,v]\in\mathcal{P}}\sigma F(B_{u},u)B_{u,v}+\mu F(B_{u},u)(v-u)\\
 & \ \ \ \ \ \ \ \ \ \ \ \ \ \ \ \ +\sigma D_{x}F(B_{u},u)\bbB_{u,v}^{\text{It\^o}}\\
 & =\lim_{|\mathcal{P}|\to0}\sum_{[u,v]\in\mathcal{P}}g(\widetilde{B}_{u})\widetilde{B}_{u,v}+Dg(\widetilde{B}_{u})\widetilde{\bbB}_{u,v}^{\text{It\^o}}+\mu F(B_{u},u)(v-u)\\
 & =\int_{s}^{t}g(\widetilde{B})\rd\widetilde{\bfB}^{\text{It\^o}}+\int_{s}^{t}\mu F(B_{u},u)du,
\end{split}
\]
which completes the proof of eqn (\ref{consistent-1}).\\
 (ii) For the second level,
\[
\begin{split} & \ \ \ \ \int_{s}^{t}f(\widetilde{B})\rrd\widetilde{\bfB}^{\text{It\^o}}\\
 & =\lim_{|\mathcal{P}|\to0}\sum_{[u,v]\in\mathcal{P}}X_{s,u}^{1}\otimes X_{u,v}^{1}+f(\widetilde{B}_{u})\otimes f(\widetilde{B}_{u})\widetilde{\bbB}_{u,v}^{\text{It\^o}}\\
 & =\lim_{|\mathcal{P}|\to0}\sum_{[u,v]\in\mathcal{P}}\left(Z_{s,u}+\int_{s}^{u}\mu F(B_{r},r)dr\right)\left(Z_{u,v}+\int_{u}^{v}\mu F(B_{r},r)dr\right)\\
 & \ \ \ \ \ \ \ \ \ \ +\sigma^{2}(F(B_{u},u))^{2}\bbB_{u,v}^{\text{It\^o}}\ \text{(Some terms go to zero here as above.)}\\
 & =\int_{s}^{t}g(\widetilde{B})\rrd\widetilde{\bfB}^{\text{It\^o}}+\int_{s}^{t}\mu Z_{s,u}^{1}F(B_{u},u)du+\int_{s}^{t}\left(\int_{s}^{u}\mu F(B_{r},r)dr\right)dZ_{u}\\
 & \ \ \ +\int_{s}^{t}\left(\int_{s}^{u}\mu F(B_{r},r)dr\right)\mu F(B_{u},u)du,
\end{split}
\]
which yields eqn (\ref{consistent-2}). 
\end{proof}

\subsection{Proof of Theorem~\ref{thm-bs-formula}}\label{appendix_thm_bs_formula}

\begin{proof} Since the price is
\[
V_{0}=e^{-rT}\mathbb{E}_{\widehat{\mathbb{P}}}\left[V_{T}\right]=e^{-rT}\mathbb{E}_{\widehat{\mathbb{P}}}\left[F(X_{T})\right],
\]
by the Girsanov theorem for fBM,
\[
\begin{split} 
\mathbb{E}_{\widehat{\mathbb{P}}}\left[F\left(X_{T}\right)\right] & =\mathbb{E}_{\widehat{\mathbb{P}}}\left[F\left(X_{0}\exp\left(\sigma B_{T}^{H}+\mu T-\frac{1}{2}\sigma^{2}T^{2H}\right)\right)\right]\\
 & =\mathbb{E}_{\widehat{\mathbb{P}}}\left[F\left(X_{0}\exp\left(\sigma\widehat{B}_{T}^{H}+rT-\frac{1}{2}\sigma^{2}T^{2H}\right)\right)\right]\\
 & =\mathbb{E}_{{\mathbb{P}}}\left[F\left(X_{0}\exp\left(\sigma{B}_{T}^{H}+rT-\frac{1}{2}\sigma^{2}T^{2H}\right)\right)\right]\\
 & =\int_{\mathbb{R}}F(X_{0}e^{\sigma T^{H}y+rT-\frac{1}{2}\sigma^{2}T^{2H}})\varphi(y)dy,
\end{split}
\]
which completes our proof.
\end{proof}


\end{document}